\documentclass[aos,MSNbibl,nameyear,dvips]{arximspdf}
\usepackage{mathbh}

\doi{10.1214/11-AOS938}
\volume{39}
\issue{6}
\pubyear{2011}
\firstpage{3262}
\lastpage{3289}

\makeatletter

\newcommand{\Nset}{{\mathbb{N}}}
\newcommand{\Rset}{\mathbb R}
\newcommand{\Cb}{\mathsf{C}_b}

\newcommand{\Uset}{\mathbb{U}}

\newcommand{\eqdef}{\stackrel{\mathrm{def}}{=}}

\newcommand{\un}{\mathbh{1}}
\newcommand{\Id}{\mathrm{I}}

\newcommand{\Xset}{\mathsf{X}}
\newcommand{\Xsigma}{\mathcal{X}}

\newcommand{\hF}{\hat{F}}

\newcommand{\tAMmean}{\mu}
\newcommand{\tAMcov}{\Gamma}

\newcommand{\PE}{\mathbb E}
\newcommand{\PP}{\mathbb P}
\newcommand{\F}{\mathcal{F}}

\newcommand{\aslim}{\stackrel{\mathrm{a.s.}}{\longrightarrow}}

\newcommand{\targetdist}{\pi}

\newtheorem{theo}{Theorem}[section]
\newtheorem{lemma}[theo]{Lemma}
\newtheorem{coro}[theo]{Corollary}
\newtheorem{prop}[theo]{Proposition}

\newproclaim{example}{Example}

\newcommand{\Xproc}{\{X_n, n \geq0\}}
\newcommand{\Yproc}{\{Y_n , n\geq0\}}

\newcommand{\rmd}{\mathrm{d}}
\newcommand{\Tset}{\Theta}

\makeatother

\begin{document}
\begin{frontmatter}

\title{Convergence of adaptive and interacting Markov chain Monte
Carlo algorithms}
\runtitle{Convergence of adaptive and interacting MCMC}

\begin{aug}
\author[A]{\fnms{G.} \snm{Fort}\corref{}\thanksref{T1}\ead[label=e1]{gersende.fort@telecom-paristech.fr}},
\author[A]{\fnms{E.} \snm{Moulines}\thanksref{T1}\ead[label=e2]{eric.moulines@telecom-paristech.fr}}
\and
\author[B]{\fnms{P.} \snm{Priouret}\ead[label=e3]{priouret@ccr.jussieu.fr}}
\runauthor{G. Fort, E. Moulines and P. Priouret}
\affiliation{LTCI, TELECOM ParisTech-CNRS,
LTCI, TELECOM ParisTech-CNRS and LPMA, Universit\'{e} Pierre et
Marie Curie}
\address[A]{G. Fort\\
E. Moulines \\
LTCI, TELECOM ParisTech-CNRS \\
46 rue Barrault \\
75634 Paris Cedex 13\\
France\\
\printead{e1}\\
\hphantom{E-mail: }\printead*{e2}}
\address[B]{P. Priouret\\
LPMA, Universit\'{e} Pierre\\
\quad et Marie Curie (P6)\\
Bo\^{i}te courrier 188 \\
75252 Paris Cedex 05\\
France\\
\printead{e3}} 
\end{aug}

\thankstext{T1}{Supported in part by the French National Research
Agency under the program ANR-08-BLAN-0218 BigMC.}

\received{\smonth{7} \syear{2010}}
\revised{\smonth{9} \syear{2011}}

%
\begin{abstract}
Adaptive and interacting Markov chain Monte Carlo
algorithms (MCMC) have been recently introduced in the literature.
These novel simulation algorithms are designed to increase the
simulation efficiency to sample complex distributions. Motivated by
some recently introduced algorithms (such as the adaptive Metropolis
algorithm and the interacting tempering algorithm), we develop a
general me\-thodological and theoretical framework to establish both the
convergence of the marginal distribution and a strong law of large
numbers. This framework weakens the conditions introduced in the
pioneering paper by Roberts and Rosenthal [\textit{J. Appl. Probab.}
\textbf{44} (2007) 458--475]. It also covers the case when the target
distribution $\pi$ is sampled by using Markov transition kernels with a
stationary distribution that differs from $\pi$.
\end{abstract}

%
\begin{keyword}[class=AMS]
\kwd[Primary ]{65C05}
\kwd{60F05}
\kwd{62L10}
\kwd{65C05}
\kwd[; secondary ]{65C40}
\kwd{60J05}
\kwd{93E35}.
\end{keyword}
\begin{keyword}
\kwd{Markov chains}
\kwd{Markov chain Monte Carlo}
\kwd{adaptive Monte Carlo}
\kwd{ergodic theorems}
\kwd{law of large numbers}
\kwd{adaptive Metropolis}
\kwd{equi-energy sampler}
\kwd{parallel tempering}
\kwd{interacting tempering}.
\end{keyword}

\end{frontmatter}

\section{Introduction}
\label{secIntroduction}
Markov chain Monte Carlo (MCMC) methods generate samples from an arbitrary
distribution $\pi$ known up to a scaling factor; see
\citet{robertcasella1999}. The algorithm consists in sampling a
Markov chain
$\Xproc$ on a general state space $\Xset$ with Markov transition
kernel~$P$
admitting $\pi$ as its unique \textit{invariant} distribution.

In most implementations of MCMC algorithms, the transition kernel $P$
of the
Markov chain depends on a tuning parameter $\theta$ defined on a space
$\Tset$
which can be either finite dimensional or infinite dimensional.

Consider, for example, the Metropolis algorithm
[\citet{metropolis1953}]. Here
$\Xset= \Rset^d$ and the stationary distribution is assumed to have a density,
also denoted by $\pi$ with respect to a measure. At the iteration $n$,
a move $Z_{n+1}=
X_n + U_{n+1}$ is proposed, where $U_{n+1}$ is drawn independently
from~$X_0,\allowbreak\ldots,X_n$ from a symmetric distribution on $\Rset^d$. This move is
accepted with probability $\alpha(X_n,Z_{n+1})$, where $\alpha(x,y)=
1 \wedge
( \pi(y)/\pi(x) )$. A frequently advocated choice of the
increment
distribution $q$ is the multivariate normal with zero-mean and covariance
matrix $(2.38^2/d)\Gamma_\star$, where $\Gamma_\star$ is the
covariance matrix
of the target distribution $\pi$ [see \citet{gelmanrobertsgilks1996}].

Of course $\Gamma_\star$ is unknown. In \citet{haariosaksmantamminen1999},
the authors have proposed an \textit{adaptive Metropolis} (AM) algorithm
in which
the covariance $\Gamma_n$ is updated at each iteration using the past
values of
the simulations [see also \citet{haariosaksman2001},
Haario et~al. (\citeyear{haariolainelehtinen2004,haariolainemirasaksman2006}),
\citet{lainetamminen2008} for applications].

The adaptive Metropolis is an example in which a parameter $\theta
_{n+1}$ is
updated at each iteration from the values of the chain $\{X_0,\ldots
,X_{n+1}\}$
and the past values of the parameters $\{\theta_0,\ldots,\theta_n\}$.
Many other
examples of such adaptive MCMC algorithms are presented in
\citet{andrieuthoms2008}, \citet{rosenthal2009} and
\citet{atchadefortmoulinespriouret2009}.

When attempting to simulate from a density with multiple modes, the Markov
kernel might mix very slowly. A useful solution to that problem is to
introduce a temperature parameter. This idea is exploited in parallel
tempering: several Metropolis algorithms are run at different
temperatures [see
\citet{geyer1991}, \citet{atchaderobertsrosenthal2010}]. One of the simulations,
corresponding to $T_1=1$ is the desired target probability
distribution. The
other simulations correspond to the family of the target distribution~$\pi^{1/T_i}$,
$i \in\{1,\ldots,K\}$, created by gradually increasing the
temperature.

The \textit{interacting tempering} algorithm, a simplified form of the
equi-energy sampler introduced \citet{kouzhouwong2006}, exploits the
parallel tempering idea. Both the algorithms run several chains in parallel,
but the interacting tempering algorithm allows more general interactions
between chains. The interacting tempering algorithm provides an example in
which the process of interest interacts with the past samples of a
family of
auxiliary processes. Other examples of such interacting schemes are presented
in \citet{andrieujasradoucetdelmoral2007} [see
also \citet{brockwelldelmoraldoucet2010}].

The two examples discussed above can be put into a common unifying framework
(see Section \ref{secMainResults}). The purpose of this work is to analyze
these general classes of adaptive and interacting MCMC. This paper complements
recent surveys on this topic by \citet{andrieuthoms2008},
\citet{rosenthal2009} and \citet{atchadefortmoulinespriouret2009}
which are
devoted to the design of these algorithms. We focus in this paper on two
problems: the ergodicity of the sampler (under which condition the marginal
distribution of the process converges to the target distribution $\pi
$) and the
strong law of large numbers (SLLN) for additive and unbounded
functionals.\vadjust{\goodbreak}

Ergodicity of the marginal distributions for adaptive MCMC
has been studied by \citet{andrieumoulines2006} for a particular
class of
samplers in which the parameter is adapted using a stochastic
approximation algorithm. These results have later been extended by
\citet{robertsrosenthal2007} to handle more general adaptation
strategies, but
under conditions which are in some respects more stringent. Most of
these works assume a form of geometric
ergodicity; these conditions are relaxed in \citet{atchadefort2008} which
addresses Markov chains with subgeometric rate of convergence.

A strong law of large number for the adaptive Metropolis algorithm was
established by \citet{haariosaksman2001} (for bounded functions and a compact
parameter space $\Theta$), using mixingales techniques; these results have
later been extended by \citet{atchaderosenthal2005} to unbounded
functions and
compact parameter space $\Theta$. The LLN for unbounded functions and
noncompact set $\Theta$ has been established recently in
\citet{saksmanvihola2010}. \citet{andrieumoulines2006} have
established the
consistency and the asymptotic normality of $n^{-1} \sum_{k=1}^n
f(X_k)$ for
bounded and unbounded functions for adaptive MCMC algorithms combined
with a
stochastic approximation procedure [see \citet{atchadefort2008} for
extensions]. The procedure involves projections on a family of increasing
compact subsets of the parameter space, and did not include the results
obtained for the AM by \citet{saksmanvihola2010}.

\citet{robertsrosenthal2007} prove a weak law of large numbers for bounded
functions for general adaptive MCMC samplers but under technical conditions
which are stringent.

The analysis of interacting MCMC algorithms started more recently and
the theory is still less developed. The original result in Kou, Zhou
and Wong [(\citeyear{kouzhouwong2006}), Theorem 2], as already noted in
the discussion paper [\citet{atchadeliu2006}, Section 3] and
carefully explained in Andrieu et al. [(\citeyear{andrieujasradoucetdelmoral2008b}),
Sec\-tion 3.1] does not amount to a proof.
\citet{andrieujasradoucetdelmoral2008b} presents a~proof of
convergence of a simple version of the interacting tempering sampler
with $K\,{=}\,2$ stages. The proofs in
\citet{andrieujasradoucetdelmoral2008b} (uniformly ergodic case)~and
in \citet{andrieujasradoucetdelmoral2008} (geometrically
ergodic case) are based on the convergence of $U$-statistics, which
explains why the results obtained for $K=2$ stages cannot easily be
extended.

SLLN was established by \citet{atchade2009} for a simple version of the
interacting tempering algorithm for a transition kernel which is geometrically
ergodic with uniformly controlled ergodicity constants, but the proof
in this
paper is not convincing [see \citet{supplement}, Section 1].

Finally, a functional Central Limit theorem
was derived in Bercu, Del~Mo\-ral and Doucet (\citeyear{bercudelmoraldoucet2009})
for a class of interacting
Markov chains for uniformly ergodic Markov kernels.

This paper aims at providing a theory weakening some of the limitations
mentioned above. Let $\{P_\theta, \theta\in\Tset\}$ be a family of
transition
kernels on $\Xset$. We address\vadjust{\goodbreak} the general framework when the target density
$\pi$ is approximated by the process $\{X_n, n \geq0\}$ such that the
conditional distribution of $X_{n+1}$ given the past is given by
$P_{\theta_n}(X_n, \cdot)$; $\{\theta_n, n\geq0\}$ is the adapted process.
There are two main contributions. First, we cover the case when the ergodicity
of the transition kernels $\{P_\theta, \theta\in\Tset\}$ is not
uniform along
the path \mbox{$\{\theta_n, n \geq0 \}$}. The second novelty is that we
address the
case when the $P_\theta$ has an invariant distribution $\pi_\theta$ depending
upon the parameter $\theta$; in this case, the adaptation has to be
such that
$\{\pi_{\theta_n}, n\geq0 \}$ converges weakly to $\pi$
(almost surely) and we
provide sufficient conditions for this property to hold based on the
(almost sure) weak convergence of the transition kernels
\mbox{$\{P_{\theta_n}, n\geq0\}$}. Such conditions are crucial in many applications where $\pi
_\theta$ is
known to exist but has no explicit expression. Therefore, to generalize the
results and include more realistic conditions, a more complex approach is
required.

The paper is organized as follows. In Section \ref{secMainResults}, we
establish the convergence of the marginal distribution and the strong
law of
large numbers for additive functionals for adaptive and interacting MCMC
algorithms. These general results are applied to a running example,
namely the
adaptive Metropolis algorithm. The novel contribution is the
application to
the convergence of the interacting tempering
algorithm [\citet{kouzhouwong2006}] in Section \ref{secEE}.\vspace*{-3pt}

\subsection*{Notation}
Let $(\Xset, \Xsigma)$ be a general state space [see, e.g.,
\citet{meyntweedie2009}, Chapter 3] and $P$
be a Markov transition kernel. $P$ acts on
bounded functions $f$ on $\Xset$ and on $\sigma$-finite positive
measures $\mu$
on $\Xsigma$ via
\[
P f(x) \eqdef\int P(x,\rmd y) f(y) ,\qquad  \mu P (A) \eqdef\int\mu
(\rmd x) P(x,A) .\vspace*{-3pt}
\]
For $n \in\Nset$, we will denote by $P^n$ the $n$-iterated transition kernel
defined by induction
\[
P^n(x,A) \eqdef\int P^{n-1}(x,\rmd y) P(y,A) = \int P(x,\rmd y)
P^{n-1}(y,A)\vspace*{-3pt}
\]
with the convention that $P^0$ is the identity kernel. For a function
$V \dvtx
\Xset\to[ 1, +\infty)$, define the $V$-norm of a function
$f\dvtx \Xset\to\Rset$ by
\[
\| f \|_{V} \eqdef\sup_{x \in\Xset} \frac{|f(x)|}{V(x)} .\vspace*{-3pt}
\]
When $V=1$, the $V$-norm is the supremum norm and will be denoted by~%
$\| f \|_{\infty}$. Let $\mathcal{L}_V$ be the set of functions such that
$\| f \|_{V} <
+\infty$. For two probability distributions $\mu_1, \mu_2$ on $\Xset
$, define
the $V$-distance
\[
\| \mu_1 -\mu_2 \|_{V} \eqdef{\sup_{\{f, \| f \|_{V} \leq1 \}}}
|\mu_1(f) - \mu_2(f)|.\vspace*{-3pt}
\]
When $V=1$, the $V$-distance is the total variation distance and is denoted
by $\|\mu_1 -\mu_2\|_{\mathrm{TV}}$.

Denote by $\Cb(\Xset)$ the
class of bounded continuous functions from $\Xset$ to~$\Rset$. Recall
that a Markov transition kernel $P$ on $(\Xset,\Xsigma)$ is (weak)
Feller if it maps $\Cb(\Xset)$ to $\Cb(\Xset)$.\vadjust{\goodbreak}

A measurable set $A \in\mathcal{A}$ on a probability space $(\Omega,
\mathcal{A}, \PP)$ is said to be a~$\PP$-full set if $\PP(A)=1$.

\section{Main results}
\label{secMainResults}
Let $(\Theta, \mathcal{T})$ be a measurable space and $(\Xset
,\Xsigma)$ a
general state space. Let $\{P_\theta, \theta\in\Theta\}$ be a
collection of Markov
transition kernels indexed by $\theta$ in~$\Theta$, which can be
either finite or
infinite dimensional. We consider a $\Xset\times\Theta$-valued
process $\{
(X_n, \theta_n), n \geq0\}$ on a filtered probability space $(\Omega,
\mathcal{A}, \{\F_n, n\geq0\}, \PP)$. It is assumed that
$(X_n,\theta_n)$ is
$\F_n$-adapted and for any bounded measurable function~$f$
%
\begin{equation}
\label{eqcondition-markov}
\mathbb{E}[f(X_{n+1}) | \F_n ] = P_{\theta_n} f(X_n) .
\end{equation}
%
%

\subsection{Ergodicity}
\label{sectheoryErgo}
For $V\dvtx \Xset\to[ 1,\infty)$ and $\theta,\theta' \in\Theta$,
denote by
$D_{V}(\theta,\theta')$ the $V$-variation of the kernels $P_{\theta
}$ and $P_{\theta'}$
%
\begin{equation}
\label{eqdefiDVnorm}
D_{V}(\theta, \theta') \eqdef\sup_{x \in\Xset} \frac{\| P_\theta
(x,\cdot) - P_{\theta'}(x,\cdot) \|_{V} }{V(x)} .
\end{equation}
When $V \equiv1$, we use the simpler notation $D(\theta, \theta')$.
Consider the following assumption:

{\renewcommand\theenumi{A\arabic{enumi}}
\renewcommand\labelenumi{\theenumi}
\begin{enumerate}
\item
\label{ALoiStationnaire}
For any $\theta\in\Theta$, there exists a probability distribution
$\pi_\theta$ such that \mbox{$\pi_\theta P_\theta= \pi_\theta$}.
\item
\label{AA2}
\begin{enumerate}[(a)]
\item[(a)]
For any $\varepsilon> 0$, there exists a nondecreasing sequence $\{
r_\varepsilon(n), n\geq0\}$ in $\Nset\setminus\{0 \}$, such that
$\limsup_{n \to\infty} r_\varepsilon(n)/n =0$ and
\[
\limsup_{n \to\infty} \PE\bigl[
\bigl\|P_{\theta_{n-r_\varepsilon(n)}}^{r_\varepsilon
(n)}\bigl(X_{n-r_\varepsilon(n)}, \cdot\bigr) - \pi_{\theta_{n-r_\varepsilon
(n)}}\bigr\|_{\mathrm{TV}} \bigr]\leq
\varepsilon.
\]
\item[(b)]
For any $\varepsilon> 0$, $\lim
_{n \to
\infty} \sum_{j=0}^{r_\varepsilon(n)-1} \PE[ D(\theta_{n
-r_\varepsilon(n) +j},
\theta_{n-r_\varepsilon(n)})] =0$, where $D$ is defined
in (\ref{eqdefiDVnorm}).
\end{enumerate}
\end{enumerate}}

Assumption \ref{AA2}(a) is implied by the \textit{containment
condition}
introduced in \citet{robertsrosenthal2007}: for any $\varepsilon\,{>}\,0$, the
sequence \mbox{$\{M_\varepsilon(X_n, \theta_n), n\,{\geq}\,0 \}$} is bounded in
probability,
where for $x \in\Xset$, $\theta\in\Theta$,
%
\begin{equation}
\label{equniform-containment-condition}
M_\varepsilon(x,\theta) \eqdef\inf\{n \geq0, \|P_\theta
^n(x,\cdot) - \pi_\theta\|_{\mathrm{TV}} \leq
\varepsilon\} .
\end{equation}
In this case, it is easily checked that \ref{AA2}(a) is
satisfied by setting $r_\varepsilon(n) = N$ for all $n \geq0$, where
$N$ is
large enough. Assumption \ref{AA2}(a) is weaker than the containment
condition, because the sequence $\{ r_\varepsilon(n), n \geq0\}$ can
grow to
infinity. This is important in applications where it is not known a
priori that
the parameter sequence $\{ \theta_n, n \geq0\}$ stays in a region
where the
ergodicity constants are controlled uniformly. Examples of such applications
are given in a toy example and a more realistic example below.

Assumption \ref{AA2}(b) requires that the amount
of change
vanishes as $n$ goes to infinity at a rate which is matched with the
rate at
which the kernel converges to stationarity. If the kernel mixes
uniformly fast
along any parameter\vadjust{\goodbreak} sequence $\{ \theta_n , n \geq0\}$, that is,
$r_\varepsilon(n)
= N$ for any $n \geq0$ for some integer $N$,~\ref{AA2}(b) is equivalent to the diminishing
adaptation condition introduced in~\citet{robertsrosenthal2007}: $\{
D(\theta_{n},
\theta_{n-1}), n\geq1\}$ converges to zero in probability at any
rate. On the
other hand, if the ergodicity is not uniform along a sequence $\{
\theta_n, n \geq
0\}$, then the rate of convergence of the adaptation should converge to zero
but with a fast enough rate. As expected, there is a~trade-off between
the rate
of convergence of the chain and the rate at which the parameter can be adapted.
This does not necessarily imply however that the parameter sequence $\{
\theta_n,
n\geq0\}$ converges to some fixed value [see, e.g.,
\citet{robertsrosenthal2007}].\vspace*{-3pt}
%
\begin{theo}
\label{theoCvgMarg}
Assume \ref{ALoiStationnaire} and \ref{AA2}. Let $f$ be a bounded
function such that $\lim_n \pi_{\theta_n}(f) = \alpha$ $\PP$-a.s.
for some constant $\alpha$. Then
\[
\lim_{n \to\infty} \PE[ f(X_n) ] = \alpha.\vspace*{-3pt}
\]
\end{theo}

The proof is in Section \ref{secprooftheoCvgMarg}. As a trivial corollary,
we have:\vspace*{-3pt}
%
\begin{coro}
\label{coroCvgMarg}
Assume \ref{ALoiStationnaire} and \ref{AA2}. Assume $\{\pi
_{\theta_n}, n\geq
0\}$ converges weakly to $\pi$ $\PP$-a.s. Then, $\lim_{n \to\infty
} \PE[
f(X_n) ] = \targetdist(f)$ for any bounded continuous function~$f$.\vspace*{-3pt}
\end{coro}

When $\pi_\theta= \pi$ for any $\theta\in\Theta$, Theorem \ref
{theoCvgMarg}
improves the results of \citet{robertsrosenthal2007} by weakening the
conditions on the transition kernels $\{P_\theta, \theta\in\Theta\}
$ (the
containment condition is not assumed to hold). The following example
shows that
ergodicity can be achieved even if the containment condition in
\citet{robertsrosenthal2007} fails, provided that the adaptation rate
is slow
enough.\vspace*{-3pt}
%
\begin{example}[(Toy example)]
Let us consider the following example introduced in
\citet{andrieumoulines2006} and thoroughly analyzed in
Andrieu and Thoms [(\citeyear{andrieuthoms2008}), Section 2]
and \citet{bairobertsrosenthal2009}. Let $\{
\theta_n,
n\geq0\}$ be a $[0,1]$-valued deterministic sequence. Consider the
nonhomogeneous Markov chain over $\Xset= \{0,1\}$ with transition matrix
%
\begin{equation}
\label{eqtransition-P-simple}
P_\theta= \left[
\matrix{
\theta& 1-\theta\cr
1-\theta& \theta}\right],\qquad
\theta\in[0,1] .
\end{equation}
For any $\theta\in[0,1]$, $\pi= [1/2,1/2]$ is a stationary
distribution; the
chain is irreducible if $\theta\in(0,1)$. In this case, for
$\varepsilon>
0$ and $\theta\in(0,1)$,
\[
M_\varepsilon(x,\theta)= \ln(\varepsilon)/{\ln}|1-2\theta| .
\]
Assume that, for $n \geq1$, $\theta_n = n^{-1/4}$. Clearly, for any
$\varepsilon> 0$, $\{ M_\varepsilon(X_n,\theta_n), n \geq0 \}$
grows to
infinity with probability $1$ and the containment condition does not
hold [see
also \citet{bairobertsrosenthal2009}, Proposition 1].

Setting $r(n)= n^{1/3}$
\[
\limsup_{n \to\infty} \PE{\bigl\|P_{\theta
_{n-r(n)}}^{r(n)}\bigl(X_{n-r(n)}, \cdot\bigr) - \pi\bigr\|_{\mathrm{TV}}}
= {\limsup_{n \to\infty}} |2 \theta_n - 1|^{r(n)}= 0\vadjust{\goodbreak}
\]
shows that A2(a) holds.
Furthermore, we have
\[
D(\theta,\theta')= {\sup_{x \in\{0,1\}}}
\|P_\theta(x,\cdot)-P_{\theta'}(x,\cdot)\|_{\mathrm{TV}} = 2
|\theta
-\theta'| .
\]
Therefore, with $\theta_n= n^{-1/4}$, $D(\theta_n,\theta_{n-1}) = O(n^{-1})$,
and A2(b) is satisfied
with $r(n)=
n^{1/3}$. Corollary \ref{coroCvgMarg} therefore applies, and the marginal
distribution converges.
\end{example}

To check \ref{AA2}(a), it is often easier to use drift
conditions. To simplify the discussion below, this paper only covers
the case
of drift inequalities for geometric ergodicity. Extensions to subgeometric
rates of convergence can be obtained following the same lines [see, e.g.,
\citet{bairobertsrosenthal2009} and \citet{atchadefort2008}] and
are left
to future work. In the geometric setting, one commonly assumes the following
simultaneous geometric drift and minorization conditions:

{\renewcommand\theenumi{A3}
\renewcommand\labelenumi{\theenumi}
\begin{enumerate}
\item
\label{Ageometric-ergodicity}
For all $\theta\in\Theta$, $P_\theta$ is $\pi$-irreducible, aperiodic
and there exist a
function $V\dvtx \Xset\to[ 1, +\infty)$, and for any $\theta\in\Theta$
there exist some constants $b_\theta< \infty$, $\delta_\theta\in
(0,1)$, $\lambda_\theta
\in(0,1)$ and a probability measure $\nu_\theta$ on $\Xset$ such that
\begin{eqnarray*}
P_\theta V &\leq&\lambda_\theta V + b_\theta, \\
P_\theta(x,\cdot) &\geq&\delta_\theta \nu_\theta(\cdot)  \un
_{\{V \leq c_\theta
\}}(x) ,\qquad  c_\theta\eqdef2 b_\theta(1-\lambda_\theta)^{-1} -1 .
\end{eqnarray*}
\end{enumerate}}

\noindent \ref{Ageometric-ergodicity} implies geometric ergodicity [see, e.g.,
\citet{meyntweedie2009}, Chapter 15]. The following proposition can be
obtained from \citet{robertsrosenthal2004}, \citet{fortmoulines2003SPA},
Douc, Moulines and Rosenthal [(\citeyear{doucmoulinesrosenthal2004}), Proposition 3]
or \citet{baxendale2005}
[see also the proof of Lemma 3 in \citet{saksmanvihola2010} for a similar
result].
%
\begin{lemma}
\label{lemBoundCandRho}
Assume \ref{Ageometric-ergodicity}. Then for any $\theta$, there
exists a
probability distribution $\pi_\theta$ such that $\pi_\theta P_\theta
= \pi_\theta$, $\pi_\theta(V)
\leq b_\theta(1-\lambda_\theta)^{-1}$ and
\[
\| P_\theta^n(x,\cdot) - \pi_\theta\|_{V} \leq C_\theta \rho
_\theta^n  V(x)
\]
for some
finite constants $C_\theta$ and $\rho_\theta\in(0,1)$. Furthermore,
there exist
positive constants $C$ and $\gamma$ such that for any $\theta\in
\Theta$,
%
\begin{equation}
\label{eqDefinitionLtheta}
L_\theta\eqdef C_\theta\vee(1-\rho_\theta)^{-1} \leq C \{
b_\theta\vee
\delta_\theta^{-1} \vee(1-\lambda_\theta)^{-1} \}^\gamma.
\end{equation}
\end{lemma}
%
%
\begin{example}[{[The adaptive Metropolis (AM) algorithm]}]
\label{secAMalgo}
We establish the ergodicity of the AM algorithm.
In this example, $\Xset= \Rset^d$ and the densities are assumed to be
w.r.t. the Lebesgue measure. For $x \in\Rset^d$, $| x |$
denotes the Euclidean norm. For $\kappa> 0$, let $\mathcal
{C}^d_\kappa$ be
the set of symmetric and positive definite $d \times d$ matrices whose minimal
eigenvalue is larger than $\kappa$. The parameter set $\Tset= \Rset
^d \times
\mathcal{C}^d_\kappa$ is endowed with the norm $| \theta|^2 \eqdef
| \mu|^2 +
\operatorname{Tr}(\Gamma^T \Gamma)$, where $\theta= ( \tAMmean,\tAMcov)$.\vadjust{\goodbreak}

At each iteration, $X_{n+1} \sim P_{\theta_n}(X_n,\cdot)$, where
$P_\theta$ is defined by
%
\begin{eqnarray}
\label{eqMetropolis-transition}
P_\theta(x,A) &\eqdef&\int_A \biggl( 1 \wedge\frac{\pi(y)}{\pi(x)}
\biggr) q_{\tAMcov}(y-x) \,\rmd y \nonumber\\[-10pt]\\[-10pt]
&&{} + \un_A(x) \biggl[ 1 - \int\biggl( 1 \wedge\frac{\pi(y)}{\pi(x)}\biggr)
q_{\tAMcov}(y-x) \,\rmd y \biggr] \nonumber\vspace*{-2pt}
\end{eqnarray}
with $q_{\tAMcov}$ the density of a
Gaussian random variable with zero mean and covariance matrix $(2.38)^2
d^{-1} \tAMcov$.
The parameter $\theta_n = (\tAMmean_n,\tAMcov_n) \in\Tset$ is the
sample mean and covariance matrix
%
\begin{eqnarray}
\label{eqAM-algorithm-mean}
\mu_{n+1} &=& \mu_{n} + \frac{1}{n+1} (X_{n+1} - \mu_{n}) ,\qquad
\mu_0 =0 , \\[-3pt]
\label{eqAM-algorithm-covariance}
\Gamma_{n+1} &=& \frac{n}{n+1}\Gamma_{n} + \frac{1}{n+1} \{
(X_{n+1} - \mu_{n})(X_{n+1} - \mu_{n})^T + \kappa\Id_d \} ,\vspace*{-2pt}
\end{eqnarray}
where $\Id_d$ is the identity matrix, $\Gamma_0 \geq0$ and $\kappa$
is a positive constant.

By construction, for any $\theta\in\Tset$, $\pi$ is the stationary
distribution
for $P_\theta$ so that \ref{ALoiStationnaire} holds with $\pi
_\theta=\pi$ for any
$\theta$. As in \citet{saksmanvihola2010}, we consider the
following assumption:
{\renewcommand\theenumi{M\arabic{enumi}}
\renewcommand\labelenumi{\theenumi}
\begin{enumerate}
\item
\label{AMPi} $\pi$ is positive, bounded, differentiable and
\[
\lim_{r \to\infty} \sup_{| x |\geq r} \frac
{x}{| x |^\rho} \cdot\nabla\log\pi(x) = - \infty\vspace*{-2pt}
\]
for some $\rho> 1$. Moreover, $\pi$ has regular contours, that is,
for some $R > 0$,
\[
\sup_{| x | \geq R} \frac{x}{| x |} \cdot\frac
{\nabla
\pi(x)}{| \nabla\pi(x) |} < 0 .\vspace*{-2pt}
\]
\end{enumerate}}

\noindent Saksman and Vihola [(\citeyear{saksmanvihola2010}), Proposition 15]
establishes a drift inequality and a minorization condition on the
kernel as in \ref {Ageometric-ergodicity}, with a drift function $V
\propto\pi^{-s}$ with $s =1/2$. Nevertheless, the generalization to an
arbitrary $s \in(0,1)$ is straightforward. Note that the function
%
\begin{equation}
\label{eqITdefinitionW}
W(x) \eqdef\pi^{-s}(x)  \| \pi^s \|_{\infty}
\end{equation}
grows faster than an
exponential under \ref{AMPi} [see, e.g., \citet{saksmanvihola2010},
Lemma 8].
Hence, Lemma \ref{lemBoundCandRho} and  Proposition 15 of
Saksman and Vihola (\citeyear{saksmanvihola2010}) both imply:\vspace*{-3pt}
%
%
\begin{lemma}
\label{lemAMgeometric-rate-of-convergence}
Assume \ref{AMPi}. For any $a \in( 0,1 ]$ and $\theta\in\Tset$,
there exist
$C_{a,\theta} < \infty$ and $\rho_{a,\theta} \in(0,1)$, such that
\[
\| P_{\theta}^k(x,\cdot) - \pi\|_{W^a} \leq C_{a,\theta}  \rho
^k_{a,\theta}  W^a(x)\qquad  \mbox{for any $x \in\Rset^d$,}\vspace*{-2pt}
\]
where $W$ is defined by (\ref{eqITdefinitionW}).
In addition, there exist finite constants $c_a,b_a$ such that
\[
C_{a,\theta} \vee(1-\rho_{a,\theta} )^{-1} \leq c_a
| \theta|^{d
\gamma/2} + b_a ,\vspace*{-2pt}
\]
where the constant $\gamma$ is defined in Lemma \ref{lemBoundCandRho}.\vspace*{-2pt}\vadjust{\goodbreak}
\end{lemma}

In Saksman and Vihola [(\citeyear{saksmanvihola2010}),
Lemma 12] it is proved that under \ref{AMPi},
the rate of growth of the parameters $\{ \theta_n, n \geq0\}$ is controlled.
Namely, for any $\tau>0$,
%
\begin{equation}
\label{eqAMStatbilitytheta}
\sup_{n \geq1} n^{-\tau} | \theta_n | <
+\infty, \qquad \PP\mbox{-a.s.}
\end{equation}
In the following lemma, we establish a control of the rate of growth of the
state of the chain $\{ X_n, n \geq0\}$.
%
\begin{lemma}
\label{lemAMStability}
Assume \ref{AMPi}. Then:

{\renewcommand\thelonglist{(\roman{longlist})}
\renewcommand\labellonglist{\thelonglist}
\begin{longlist}
\item\label{propAMBoundW} $\PE[W(X_n) ] \leq\PE
[W(X_0)
] + n b $.
\item\label{propAMStabilityW} For any $t >0$ and any $\tau>0$, there
exists a constant $C_{t,\tau}$ such that for any $n \geq0$,
\[
\PE\bigl[W(X_n) \un_{\sup_{k \leq n-1} k^{-\tau} | \theta_k |
\leq t }\bigr] \leq
\PE[W(X_0) ] + C_{t,\tau}  n^{\tau d \gamma/2} ,
\]
where $\gamma$ is defined in Lemma \ref{lemBoundCandRho}.
\item\label{propAMStabilityW-ps} If $\PE[W(X_0) ]\,{<}\,
{+}\infty$,
for any $\tau\,{>}\,0$, $\sup_{n \geq1} n^{-1-\tau} W(X_n)\,{<}\,{+}\infty$,~\mbox{$\PP$-a.s.}
\end{longlist}}
\end{lemma}

The proof of this lemma is given in Section \ref{secProofAM}. By combining
Lemma \ref{lemAMgeometric-rate-of-convergence} and
Lemma \ref{lemAMStability}, we prove \ref{AA2}(a): as a consequence of Lemma \ref
{lemAMgeometric-rate-of-convergence}, it
holds for any $\tau>0$ such that $r > \tau d \gamma/2$ and for any $t>0$
%
\begin{equation}
\label{eqpropAMUniformBound-1}
\limsup_{n \to\infty}   \sup_{\theta\in\Tset, | \theta| \leq
t n^\tau}
\sup_{x \in\Rset^d, W(x) \leq t n^{1+\tau}} \bigl\|P_\theta
^{\lfloor n^{r} \rfloor} (x,\cdot) - \pi\bigr\|_{\mathrm{TV}} = 0 ,
\end{equation}
where $\lfloor\cdot\rfloor$ denotes the lower integer part.
For $t > 0$, set
\[
\Omega_t \eqdef\Bigl\{\omega\dvtx \sup_{n \geq1}
n^{-\tau} | \theta_n | \leq t , \sup_{n \geq1} n^{-1 -
\tau}
W(X_n) \leq t \Bigr\}.
\]
Equation (\ref{eqAMStatbilitytheta}) and
Lemma \ref{lemAMStability}\ref{propAMStabilityW-ps} show that
$\lim_{t
\to\infty} \PP( \Omega_t ) = 1$. Set $r(n)= \lfloor n^{r} \rfloor
$. The
Fatou lemma and the monotone convergence theorem show that
\begin{eqnarray*}
&&
\limsup_{n \to\infty} \PE\bigl[
\bigl\|P^{r(n)}_{\theta_{n-r(n)}}\bigl(X_{n-r(n)},\cdot\bigr) - \pi\bigr\|_{\mathrm
{TV}}
\bigr] \\
&&\qquad
\leq\PE\Bigl[ \limsup_{n \to\infty} \bigl\|P^{r(n)}_{\theta
_{n-r(n)}}\bigl(X_{n-r(n)},\cdot\bigr) - \pi\bigr\|_{\mathrm{TV}} \Bigr] \\
&&\qquad\leq\lim_{t \to\infty} \PE\Bigl[ \limsup_{n \to\infty}
\bigl\|P^{r(n)}_{\theta_{n-r(n)}}\bigl(X_{n-r(n)},\cdot\bigr) -
\pi\bigr\|_{\mathrm{TV}} \un
_{\Omega_t}
\Bigr] = 0 .
\end{eqnarray*}
%
Therefore, \ref{AA2}(a) is satisfied whereas
clearly the
uniform containment condition [see
(\ref{equniform-containment-condition})]
seems to be very challenging to check.

Consider now \ref{AA2}(b). It is proved in
Andrieu and Moulines
[(\citeyear{andrieumoulines2006}), Lem\-ma~13] that for any $(\theta,\tilde
{\theta}) \in
\Theta^2$ and $a \in[0,1]$, $D_{W^a}(\theta,\tilde{\theta}) \leq2
d \kappa^{-1}
| \tAMcov- \tilde{\Gamma} |$. By definition of $\tAMcov
_n$ [see
(\ref{eqAM-algorithm-covariance})], we have for any $m < n$,
\begin{eqnarray*}
D_{W^a}(\theta_n,\theta_{n-m}) &\leq&\frac{2d\kappa^{-1}}{n}
\Biggl(2 \kappa m d +
\frac{m}{n-m} \sum_{j=0}^{n-m-1} | X_{j+1}-\tAMmean_j |^2
\\[-3pt]
&&\hspace*{97pt}{} +
\sum_{j=n-m}^{n-1} | X_{j+1} -\tAMmean_{j} |^2 \Biggr) .
\end{eqnarray*}
By definition of the empirical mean $\tAMmean_k$ [see
(\ref{eqAM-algorithm-mean})] there exists a constant $C'$ such that
$|\tAMmean_k| \leq C'  \{k^{-1} \sum_{j=1}^k |X_j|^2\}^{1/2}$; under
\ref{AMPi}, $\liminf_{| x | \to\infty} \ln W(x)/ | x | >
0$\vspace*{1pt} [see the proof of Lemma 8 in
\citet{saksmanvihola2010}]. Therefore, there exists a constant $C$
such that
%
\begin{eqnarray}
\label{eqAMControleDw}
&&
D_{W^a}(\theta_n,\theta_{n-m})\nonumber \\[-2pt]
&&\qquad
\leq C  \frac{m}{n} \Biggl\{ 1 + \frac{(1 + \ln(n-m) )}{n-m}
\sum_{j=1}^{n-m} \ln^2W(X_j)\\[-2pt]
&&\hspace*{74pt}{} + \frac{(1 + \ln(n) )}{m}
\sum_{j=n-m}^{n} \ln^2 W(X_j)\Biggr\}.
\nonumber
\end{eqnarray}
The proof of \ref{AA2}(b) now relies on the control of moments for the
r.v.\break $\{ \ln^2 W(X_j)$, $j \geq0\}$. Lemma
\ref{lemAMStability}\ref{propAMBoundW} and Jensen's inequality show
that the moment $\PE[ \ln^2 W(X_n)]$ increases at most as $\ln^2
n$. Then there exists a constant $C$ such that for any $m \leq n$ and
for any $a \in[0,1]$,
\[
\PE[D_{W^a}(\theta_n,\theta_{n-m}) ] \leq C  m  \frac
{\ln^{3}(n)}{n}
\PE[W(X_0) ] .
\]

Then,\vspace*{1pt} for any $r \in(0,1/2)$, $\lim_{n \to+\infty} \sum
_{j=0}^{\lfloor
n^{r}\rfloor-1}  \PE[ D(\theta_{n-\lfloor n^{r}\rfloor+j},
\theta_{n-\lfloor
n^{r}\rfloor})] = 0$ and \ref{AA2}(b)
holds. Combining the results above yields:
\end{example}
%
\begin{theo}
\label{theoAMErgodicity}
Assume \ref{AMPi} and $\PE[W(X_0) ] < +\infty$. Then,
for any
boun\-ded function $f$, $\lim_{n \to\infty} \PE[f(X_n)] = \pi(f)$.
\end{theo}

\subsection{Strong law of large numbers for additive functionals}
\label{sectheoryLGN}

In this section, a~strong law of large numbers (SLLN) is established.
The main result of this section is Theorem \ref{theoLGNgal} which
provides a SLLN for a special class of additive functionals. To that
goal, \ref{Ageometric-ergodicity} is assumed to hold (which implies
\ref{ALoiStationnaire}, see Lemma \ref{lemBoundCandRho}), and it is
required to strengthen the diminishing adaptation and the stability
conditions.

{\renewcommand\theenumi{A4}
\renewcommand\labelenumi{\theenumi}
\begin{enumerate}
\item
\label{AStrongContainmentLGN}
$ \sum_{k =1}^{\infty} k^{-1} ( L_{\theta_k} \vee L_{\theta
_{k-1}} )^6
D_{V}(\theta_k, \theta_{k-1})  V(X_k) < +\infty$ $\PP$-a.s., where
$D_V$ and
$L_\theta$ are defined in (\ref{eqdefiDVnorm}) and
(\ref{eqDefinitionLtheta}).
\end{enumerate}}

{\renewcommand\theenumi{A5}
\renewcommand\labelenumi{\theenumi}
\begin{enumerate}
\item
\label{AUniformBoundLGN}
\begin{enumerate}[(a)]
\item[(a)] 
$\limsup_n \pi_{\theta_n}(V)
< +\infty$,
$\PP$-a.s.
\item[(b)] 
For some $\alpha> 1$,
$\sum_{k=0}^\infty(k+1)^{-\alpha}  L_{\theta_{k}}^{2\alpha}
P_{\theta_{k}}V^{\alpha}(X_{k}) < + \infty$, $\PP$-a.s.\vadjust{\goodbreak}
\end{enumerate}
\end{enumerate}}

\noindent Here again, these conditions balance the rate at which the transition
kernel~$P_\theta$ converges to stationarity and the adaptation speed. This is
reflected in
the condition~\ref{AStrongContainmentLGN}: $(L_{\theta_k} \vee
L_{\theta_{k-1}})$ is
related to the rate of convergence of the kernels~$P_{\theta_k}$ and~%
$P_{\theta_{k-1}}$
to stationarity and $D_V(\theta_k, \theta_{k-1})$ reflects the
adaptation speed.

\begin{theo}
\label{theoLGNgal}
Assume \ref{Ageometric-ergodicity}, \ref{AStrongContainmentLGN} and
\ref{AUniformBoundLGN}. Let $F\dvtx \Xset\times\Theta\to\Rset$ be a
measurable function such that:
{\renewcommand\thelonglist{(\roman{longlist})}
\renewcommand\labellonglist{\thelonglist}
\begin{longlist}
\item\label{theoLGNgalitem1} $\sup_\theta\| F(\cdot,\theta) \|
_{V} <
+\infty$,
\item\label{theoLGNgalitem2} $ \sum_{k=1}^{\infty} k^{-1}
L_{\theta_{k-1}}^2
\| F(\cdot, \theta_k) - F(\cdot, \theta_{k-1}) \|_{V}  V(X_k)<
+\infty$ $\PP$-a.s.,
\item\label{theoLGNgalitem3} $\lim_{n \to\infty} \int\pi
_{\theta_n}(\rmd x)
F(x,\theta_n)$ exists $\PP$-a.s.
\end{longlist}}

\noindent Then,
\[
\lim_{n \to\infty} \frac{1}{n} \sum_{k=0}^{n-1} F(X_k, \theta_k)
= \lim_{n \to
\infty} \int\pi_{\theta_n}(\rmd x) F(x,\theta_n) ,\qquad
\PP \mbox{-a.s.}
\]
\end{theo}

The proof is in Section \ref{secprooftheoLGNgal}. When the function
$F$ does
not depend upon~$\theta$, this theorem becomes
the following.
%
\begin{coro}
\label{coroLNGgal}
Assume \ref{Ageometric-ergodicity}, \ref{AStrongContainmentLGN} and
\ref{AUniformBoundLGN}. Let $f\dvtx \Xset\to\Rset$ be a measurable function
such that $\| f \|_{V} < +\infty$ and $\lim_{n \to\infty} \pi
_{\theta_n}(f)$
exists $\PP$-a.s. Then, $n^{-1} \sum_{k=0}^{n-1} f(X_k)
\stackrel{\mathit{a.s.}}{\longrightarrow}\lim_n
\pi_{\theta_n}(f)$.
\end{coro}
%
\begin{example}[(Toy example: law of large numbers)]
For $\theta\in(0,1)$, the constants $C_\theta$ and $\rho_\theta$ (see
Lemma \ref{lemBoundCandRho}) are, respectively, equal to $1$ and $ |1-2
\theta |$ and $V=1$. This implies that $ L_\theta= 1/(2\theta)$ if
$\theta\leq 1/2$ and $1/(2(1-\theta))$ otherwise. Therefore A3 is
satisfied since $ \sum_{k =1}^{\infty} k^{-1} \theta_{k}^{-3}
|\theta_{k-1} - \theta_k| < +\infty$ when $\theta_k = k^{-1/4}$.
Assumption A4(a) is automatically satisfied because the stationary
distribution does not depend on $\theta$. Assumption A4(b) is satisfied
for any $\alpha> 4/3$ because in such case $\sum _{k=1}^\infty (k^{-1}
\theta_k)^\alpha< \infty$. By Theorem \ref{theoLGNgal}, the SLLN is
satisfied for this nonhomogeneous Markov chain.
\end{example}

The stated assumptions are very general and, when applied to some
specific settings, can be simplified. For example, in many interesting
examples (see, e.g., Section \ref{secEE}), it is known that $\limsup_{n
\to\infty} L_{\theta_n} < \infty$, $\PP$-a.s. and for some $\alpha> 1$,
$\sup _{n \geq0} \PE[ V^\alpha(X_n)] < \infty$. Under\vspace*{1pt}
these assumptions, it is straightforward to establish the following
corollary:
%
\begin{coro}
\label{coroLNGgalEE}
Assume \ref{Ageometric-ergodicity} and:
{\renewcommand\thelonglist{(\roman{longlist})}
\renewcommand\labellonglist{\thelonglist}
\begin{longlist}
\item\label{coroLNGgalEEitem1} $\limsup_{n \to\infty} L_{\theta
_n} < \infty$
and $\limsup_{n \to\infty} \pi_{\theta_n}(V) < +\infty$, $\PP
$-a.s.,
\item\label{coroLNGgalEEitem2} there exists $\alpha>1$ such that
$\sup_{k
\geq0} \PE[V^\alpha(X_k)] < +\infty$,
\item\label{coroLNGgalEEitem3} $\sum_{k=1}^\infty k^{-1}
D_{V}(\theta_k, \theta_{k-1})  V(X_k) < +\infty$ $\PP$-a.s.\vadjust{\goodbreak}
\end{longlist}}

\noindent Let $f\dvtx \Xset\to\Rset$ be a measurable function such that $ \| f \|_{V}<
+\infty$ and\break $\lim_{n \to\infty} \pi_{\theta_n}(f)$ exists $\PP
$-a.s. Then,
$n^{-1} \sum_{k=0}^{n-1} f(X_k) \stackrel{\mathit{a.s.}}{\longrightarrow}\lim_{n \to\infty}\pi
_{\theta_n}(f)$.
\end{coro}
%
\begin{example}[(AM: law of large numbers)]
Application of the above criteria yields the SLLN for the AM algorithm. This
result has recently been obtained by \citet{saksmanvihola2010}.

Let $a \in(0,1)$ and set $W(x) \eqdef\pi^{-s}(x) \| \pi^{s} \|
_{\infty}$ for $s
\in(0,1)$. We prove that a (strong) LLN holds for any function $f$ in
$L_{W^a}$. We choose $\tau>0$ small enough so that
%
\begin{equation}
\label{eqContraintesVarepsilon}
(1-a) > \tau(a + 3 d \gamma) ,\qquad  1/a-1 > \tau d \gamma
(1/a+1/2) ,
\end{equation}
where $\gamma$ is given by Lemma \ref{lemBoundCandRho}. Consider
\ref{AStrongContainmentLGN}. By
Lemma \ref{lemAMgeometric-rate-of-convergence} and
(\ref{eqAMStatbilitytheta}), there exists a r.v. $U_1$, $\PP$-a.s.
finite
such that $L_{\theta_k} \vee L_{\theta_{k-1}} \leq U_1  k^{\tau d
\gamma/2}$. By
(\ref{eqAMControleDw}) and
Lem\-ma~\ref{lemAMStability}\ref{propAMStabilityW-ps}, there
exists a
r.v. $U_2$, $\PP$-a.s. finite such that $D_{W^a}(\theta_k,\theta
_{k-1}) \leq U_2
k^{-1} \ln^3 k$. Finally, applying
Lemma \ref{lemAMStability}\ref{propAMStabilityW-ps} again,
there exists
a r.v. $U_3$, $\PP$-a.s. finite such that $W^a(X_k) \leq U_3
k^{a(1+\tau)}$.
Combining these inequalities show that there exists a r.v. $U$, $\PP
$-a.s.
finite such that
\[
\sum_k k^{-1} ( L_{\theta_k} \vee L_{\theta_{k-1}} )^6
D_{W^a} (\theta_k,
\theta_{k-1})  W^a(X_k) \leq U \sum_k k^{2 - a - \tau(a + 3 d
\gamma)} \ln^3 k
,
\]
thus showing \ref{AStrongContainmentLGN} [observe that the RHS is
finite by
definition of $\tau$, equation~(\ref{eqContraintesVarepsilon})]. The
proof of
\ref{AUniformBoundLGN}(b) could rely on
the same
inequalities in the case $a \in(0,1/2)$. Nevertheless, a SLLN can be
established for larger values of $a$ by using the bound on $W(X_n)$
given by
Lemma \ref{lemAMStability}\ref{propAMStabilityW} which
improves on
Lemma~\ref{lemAMStability}\ref{propAMStabilityW-ps}. Set
$\Omega_t
\eqdef\{ \sup_{n \geq1} n^{-\tau}|\theta_n| \leq t \}$. By
Lemma \ref{lemAMStability}, $\lim_{t \to+\infty} \PP(\Omega
_{t})\,{\uparrow}\,1$
and \ref{AUniformBoundLGN}(b) holds provided
$\sum_{k \geq1} k^{-1/a} L_{\theta_{k-1}}^{2/a} P_{\theta_{k-1}} W(X_k)
\un_{\Omega_{t}}$ is finite $\PP$-a.s. for any $t>0$.
Lem\-mas~\ref{lemAMgeometric-rate-of-convergence} and
\ref{lemAMStability}\ref{propAMStabilityW} imply that there
exists a constant $C_{t}$ such that
\[
\PE\biggl[ \sum_k k^{-1/a} L_{\theta_{k-1}}^{2/a} P_{\theta_{k-1}} W(X_k)
\un_{\Omega_{t}} \biggr] \leq C_{t}  \sum_k k^{-1/a + \tau d
\gamma
(1/a+1/2) } .
\]
The RHS is finite by definition of $\tau$ [see
(\ref{eqContraintesVarepsilon})].
\end{example}

The above discussion is summarized in the following theorem.
%
\begin{theo}
\label{theoAMLLN} Assume \ref{AMPi} and $\PE[W(X_0)
] < +\infty$. Then, for any $a
\in(0,1)$ and any function $f \in\mathcal{L}_{W^a}$, $n^{-1} \sum
_{k=1}^n f(X_k)
\stackrel{\mathit{a.s.}}{\longrightarrow}\pi(f)$.
\end{theo}
%

\subsection{Almost sure convergence of the invariant distributions}
\label{secaspi}
When the stationary distribution $\pi_{\theta}$ is not explicitly known,
convergence of the sequence $\{\pi_{\theta_n}, n\geq0\}$ has to be
obtained from
the convergence of the transition kernels $\{P_{\theta_n}, n\geq0\}$.
We propose
below a set of sufficient conditions allowing to prove the almost sure
convergence of $\{\pi_{\theta_n}(f), n\geq0\}$ for continuous
functions $f$. The
proof of Theorem \ref{theoaspi} is in Section \ref{secproofaspi}.\vadjust{\goodbreak}
%
\begin{theo}
\label{theoaspi}
Assume that $\Xset$ is a Polish space. Assume \ref
{Ageometric-ergodicity} and:
{\renewcommand\thelonglist{(\roman{longlist})}
\renewcommand\labellonglist{\thelonglist}
\begin{longlist}
\item\label{theoaspiitem1} $ \limsup_{n \to\infty} L_{\theta_n}
< \infty$ $\PP$-a.s. where $L_\theta$ is given by (\ref
{eqDefinitionLtheta}),
\item\label{theoaspiitem3} for any function $f$ in $\Cb(\Xset)$,
the class
of functions $\{P_{\theta}f, \theta\in\Theta\}$ is equicontinuous,
\item\label{theoaspiitem2} there
exists $\theta_\star\in\Theta$ and for any $x \in\Xset$, a $\PP
$-full set $\Omega_x$ such that for any $\omega\in\Omega_x$, $\{
P_{\theta_n(\omega)}(x,\cdot),
n\geq0\}$ converges weakly to $P_{\theta_\star}(x,\cdot)$.
\end{longlist}}

\noindent Then, there exists a $\PP$-full set $\Omega_0$ such that, for any any
$\omega
\in\Omega_0$ and $f \in\Cb(\Xset)$, $\pi_{\theta_n(\omega)}(f)
\stackrel{\mathit{a.s.}}{\longrightarrow}
\pi_{\theta_\star}(f)$ (or, equivalently, for any $\omega\in
\Omega_0$, $\pi_{\theta_n(\omega)}$ converges weakly to $\pi
_{\theta_\star}$).
\end{theo}

Note that the weak convergence implies that for any $\omega\in\Omega
_0$ and
for any set $A$ such that $\pi_{\theta_\star}(\partial A)=0$ where
$\partial A$
denotes the boundary of $A$, $\lim_n \pi_{\theta_n(\omega)}(A) =
\pi_{\theta_\star}(A)$.

Theorem \ref{theoaspi} might be
seen as an extension of the classical results on the continuity of the
perturbations of the spectrum and eigenprojections; but it is stated under
assumptions that are weaker than what is usually assumed
[\citet{kato1980}, Theorem 3.16]. The difference stems from the fact that condition
\ref{theoaspiitem2} does not imply the convergence of $P_\theta$ to
$P_{\theta_\star}$ in operator norm. This is crucial to deal with the
interacting
tempering algorithm (see Section \ref{secEE}).

Condition \ref{theoaspiitem2} of Theorem \ref{theoaspi} is
certainly the
most difficult to check. In the case, it is known that for any function
$f \in
\Cb(\Xset)$, there exists a $\PP$-full set $\Omega_{x,f}$ such that
for any
$\omega\in\Omega_{x,f}$, $ \lim_nP_{\theta_n(\omega)}f(x) =
P_{\theta_\star}f(x)$,
then the existence of a $\PP$-full set, uniform in $f$ for $f \in\Cb
(\Xset)$,
relies on the characterization of the weak convergence by a separable
class of
functions [see \citet{dudley2002}, Theorem 11.4.1, and
Proposition \ref{propITaspi} below for an example].

\section{Convergence of the interacting tempering (IT) algorithm}
\label{secEE}
We consider the interacting tempering algorithm, which is a simplified
form of
the equi-energy sampler by \citet{kouzhouwong2006}.

Assume that $\Xset$ is a Polish space equipped with its Borel $\sigma$-field
$\Xsigma$. Let~$\pi$ be the target density w.r.t. a measure $\mu$ on
$(\Xset,
\Xsigma)$. Denote by $K$ the number of different temperature levels,
$T_1=1 <
T_2 < \cdots< T_K$. For $k \in\{1,\ldots, K-1\}$, let $P^{(k)}$ be a
transition kernel on $(\Xset, \Xsigma)$ with unique invariant
distribution~$\pi^{1/T_k}$. Fix $\upsilon\in(0,1)$ the probability of interaction.

We denote by $X^{(k)} = (X_n^{(k)})_n$ the sampled values at each temperatures~$T_k$. The chains are defined by induction on $k$: given the past of the
process~$X^{(k+1)}$ up to time $n$, and the current value $X_{n}^{(k)}$
of the
current process $X^{(k)}$, we define $X_{n+1}^{(k)}$ as follows:
\begin{enumerate}
\item with probability $(1-\upsilon)$, the state $X^{(k)}_{n+1}$ is sampled
using the Markov kernel $P^{(k)}(X_n^{(k)},\cdot)$,
\item with probability $\upsilon$, a tentative state $Z_{n+1}$ is
drawn at random from
the past $\{X^{(k+1)}_\ell, \ell\leq n\}$. This move is accepted
with probability
$1 \wedge(\pi(X_n^{(k)})/\pi(Z_{n+1}))^{T^{-1}_{k+1}-T_k^{-1}}$.
\end{enumerate}
We consider first the case $K=2$. We will then address the general case (see
Theorem \ref{theoEEKstage} below). For notational simplicity, we set
$T_2= T
> 1$ and $P^{(1)}=P$. Denote by $\Tset$ the set of the probability
measures on
$(\Xset, \Xsigma)$. For any distribution $\theta\in\Tset$, define
the transition
kernel $P_\theta(x,\cdot) \eqdef(1-\upsilon) P(x,\cdot) + \upsilon
K_\theta(x,\cdot)$,
where, for any $A \in\Xsigma$,
%
\begin{equation}
\label{eqdefinition-Kt}
K_\theta(x,A) \eqdef\int_A \alpha(x,y) \theta(\rmd y) + \un_A(x)
\int\{ 1- \alpha(x,y) \} \theta(\rmd y)
\end{equation}
with
%
\begin{equation}
\label{eqacceptance-ratio}
\alpha(x,y)= 1 \wedge\frac{\pi(y) \pi^{1/T}(x)}{ \pi(x) \pi
(y)^{1/T}} = 1 \wedge\frac{\pi^\beta(y)}{\pi^\beta(x)} ,\qquad
\beta\eqdef1 - \frac{1}{T} \in(0,1) .
\end{equation}
Denote by $\Yproc$ the process run at the temperature $T$. It is not assumed
that $\Yproc$ is a Markov chain. We simply assume that, for any bounded
continuous function $f$, $n^{-1} \sum_{k=1}^n f(Y_k) \to\theta_\star
(f)$ a.s. where
$\theta_\star$ is the probability distribution on $(\Xset,\Xsigma)$
with density
(w.r.t. $\mu$) proportional to $\pi^{1/T}$. We consider the process
$\Xproc$
defined, for each $n \geq0$ and any bounded function $f\dvtx \Xset\to
\Rset$,
\[
\mathbb{E}[f(X_{n+1}) | \F_n ]= P_{\theta_n} f(X_n) \qquad  \mbox{where }
\theta_n(f)= (n+1)^{-1} \sum_{k=0}^n f(Y_k) .
\]
Since, by construction, $\pi P_{\theta_\star} = \pi$, it is expected
that the
marginal distribution of $X_k$ as $k$ goes to infinity converges to
$\pi$. To
go further, some additional assumptions are required:
{\renewcommand\theenumi{I\arabic{enumi}}
\renewcommand\labelenumi{\theenumi}
\begin{enumerate}
\item\label{EESPi}
\label{EBoundPi} $\pi$ is a continuous positive density on $\Xset$
and $\| \pi\|_{\infty} < +\infty$.
\item
\label{EKernelP}
\begin{enumerate}[(a)]
\item[(a)]
$P$ is a $\pi$-irreducible aperiodic Feller
transition kernel on
$(\Xset, \Xsigma)$ such that $\pi P = \pi$.
\item[(b)] 
There exist $\tau\in(0,1/T)$, $\lambda\in
(0,1)$ and
$b<+\infty$ such that
%
\begin{equation}
\label{eqdefinitionW}
P W \leq\lambda W + b  \qquad\mbox{with }
W(x) \eqdef\bigl( \pi(x) / \| \pi\|_{\infty} \bigr)^{-\tau} .
\end{equation}
\item[(c)]
For any $p \in(0,\| \pi\|_{\infty})$, the
sets $\{
\pi\geq p \}$ are
$1$-small (w.r.t. the transition kernel $P$).
\end{enumerate}
\end{enumerate}}

\noindent When $\Xset\subseteq\Rset^d$ and $P$ is a symmetric random-walk
Metropolis (SRWM) algorithm then $\pi P = \pi$ and $P$ is
$\pi$-irreducible [\citet{mengersentweedie1996}, Lem\-ma~1.1]. If
in addition the proposal density is continuous on $\Xset$ then, since
$\pi$ is positive and continuous on $\Xset$, any compact set of $\Xset$
is $1$-small [\citet{mengersentweedie1996}, Lemma 1.2]. Therefore,
the transition kernel of a SRWM algorithm satisfies \ref{EKernelP}(a)
and \ref{EKernelP}(c).\vadjust{\goodbreak}

Drift conditions of the form \ref{EKernelP}(b) for the SRWM
algorithm on $\Xset\subseteq\Rset^d$ are discussed in
\citet{robertstweedie1996}, \citet{jarnerhansen2000} and
\citet{saksmanvihola2010}. Under conditions which imply that the target
density $\pi$ is super-exponential and have regular contours (see
\ref{AMPi}), \citet{jarnerhansen2000} and \citet
{saksmanvihola2010} show
that any functions proportional to $\pi^{-s}$ with $s \in(0,1)$
satisfies a
Foster--Lyapunov drift inequality [\citet{jarnerhansen2000},
Theorems 4.1 and 4.3]. Under this condition,~\ref{EKernelP}(b)
is satisfied with any $\tau$ in the interval $(0, 1/T)$.

Stability conditions on the auxiliary process $\{Y_n, n\geq0\}$ are also
required.

{\renewcommand\theenumi{I3}
\renewcommand\labelenumi{\theenumi}
\begin{enumerate}
\item
\label{EProcY}
\begin{enumerate}[(a)]
\item[(a)] 
$\theta_\star(W) < +\infty$ and for any
continuous
function $f$ in $\mathcal{L}_{W}$, $\theta_n(f) \aslim\theta_\star(f)$.
\item[(b)] 
$\sup_n \PE[W(Y_n) ]< + \infty$.
\end{enumerate}
\end{enumerate}}

The following proposition is the key-ingredient to prove the
convergence of the
IT sampler. Under the stated assumptions, we prove that the transition kernels
$\{P_\theta, \theta\in\Theta\}$ satisfy a Foster--Lyapunov drift
inequality and a
minorization condition. The proof of Proposition \ref{propkey-iMCMC} is
adapted from Atchad{\'e} [(\citeyear{atchade2009}), Lemma 4.1]; a detailed proof is given
in \citet{supplement}, Section 2.
%
\begin{prop}
\label{propkey-iMCMC}
Assume \ref{EESPi} and \ref{EKernelP}. Then, there exist $\tilde
\lambda
\in(0,1)$, \mbox{$\tilde b < \infty$}, such that, for any $\theta\in\Theta$,
%
\begin{equation}
\label{eqEEdriftineq}
P_\theta W(x) \leq\tilde\lambda W (x)+ \tilde b \theta(W) .
\end{equation}
In addition, for any $p \in(0,\| \pi\|_{\infty})$, the level sets
$\{\pi
\geq p
\}$ are $1$-small w.r.t. the transition kernels $P_\theta$ and the minorization
constant does not depend upon $\theta$.
\end{prop}
%
\begin{coro}
\label{corokey-iMCMC}
Assume \ref{EESPi}, \ref{EKernelP}, \ref{EProcY} and $\PE
[W(X_0)]< +
\infty$. Then:
{\renewcommand\thelonglist{(\roman{longlist})}
\renewcommand\labellonglist{\thelonglist}
\begin{longlist}
\item\label{corokey-iMCMCitem1} $\sup_{n \geq0} \PE[W(X_n)
] <
+\infty,$
\item\label{corokey-iMCMCitem2} $\limsup_{n \to\infty} L_{\theta
_n} < +\infty$
$\PP$-a.s., where $L_\theta$ is defined by (\ref{eqDefinitionLtheta}).
\end{longlist}}
\end{coro}

The proof of Corollary \ref{corokey-iMCMC} is in
Section \ref{secproofcorokey-iMCMC}. As a consequence of
Proposition~\ref{propkey-iMCMC}, the transition kernel $P_\theta$
possesses an
(unique) invariant distribution $\pi_\theta$. Ergodicity and SLLN for additive
functionals both require the a.s. convergence of $\pi_{\theta_n}(f)$
(see
Theorems \ref{theoCvgMarg} and \ref{theoLGNgal}). Nevertheless, in this
example, $\pi_\theta$ does not have an explicit expression. The proof
of the
following proposition is postponed in Section \ref{secproofpropITaspi}.
%
\begin{prop}
\label{propITaspi}
Assume \ref{EESPi}, \ref{EKernelP}, \ref{EProcY} and $\PE
[W(X_0)]< +
\infty$. Then, the conditions of Theorem \ref{theoaspi} hold and for any
bounded continuous function~$f$, $\lim_n \pi_{\theta_n}(f) = \pi
(f)$ $\PP$-a.s.
\end{prop}

We now address the convergence of the marginals.\vadjust{\goodbreak}
%
\begin{theo}
\label{theoEECvgMarginal}
Assume \ref{EESPi}, \ref{EKernelP}, \ref{EProcY} and $\PE
[W(X_0)]< +
\infty$. Then, for any bounded continuous function $f$, $ \lim_n
\PE[f(X_n) ] = \pi(f)$.
\end{theo}
\begin{pf}
We check the assumptions of Corollary \ref{coroCvgMarg}. By
Corollary \ref{corokey-iMCMC}\ref{corokey-iMCMCitem1}, $\{
W(X_n), n\geq
0\}$ is bounded in probability. Furthermore,
Corollary \ref{corokey-iMCMC}\ref{corokey-iMCMCitem2} implies that
$\limsup_n C_{\theta_n}< +\infty$ $\PP$-a.s. and $\limsup_n \rho
_{\theta_n}
<1$ $\PP$-a.s. This pro\-ves~\ref{AA2}(a).

The next step is to establish \ref{AA2}(b).
Since, for any bounded function $f$, $\theta_{n+m}(f) = (n+m+1)^{-1}
\sum_{k=n+1}^{n+m} f(Y_k) + (n+1) (n+m+1)^{-1} \theta_{n}(f)$, we have
\[
| P_{\theta_{n+m}} f(x) - P_{\theta_{n}} f(x) | \leq
{\sup_{y,z \in\Xset} }|
f(y) -f(z) |  \|\theta_{n+m} - \theta_{n}\|_{\mathrm{TV}} \leq
\frac{2
\| f \|_{\infty} m}{n+m+1}
.
\]
Consequently, $D(\theta_{n+m},\theta_{n})$ is deterministically
bounded by a sequence
converging to zero. We have
\[
\sum_{j=0}^{r_\varepsilon(n) -1} \PE\bigl[D\bigl(\theta
_{n-r_\varepsilon(n) +j},
\theta_{n-r_\varepsilon(n)}\bigr) \bigr] \leq2
\frac{r_\varepsilon^2(n)}{n-r_\varepsilon(n)}
\]
thus proving \ref{AA2}(b) with any
sequence of the
form $r_\varepsilon(n) = n^r$ with $r <1/2$.

Finally, Proposition \ref{propITaspi} proves the convergence of
$\pi_{\theta_n}(f)$ for any bounded continuous function $f$.
\end{pf}

We now state the strong law of large numbers for the IT sampler.
%
\begin{theo}
\label{theoEECvgResults}
Assume \ref{EESPi}, \ref{EKernelP}, \ref{EProcY} and $\PE
[W(X_0)]< +
\infty$. Then:
{\renewcommand\thelonglist{(\roman{longlist})}
\renewcommand\labellonglist{\thelonglist}
\begin{longlist}
\item for any measurable set $A$ such that $\int_{\partial A} \pi
\,\rmd\mu= 0$
where $\partial A$ is the boundary of $A$,
\[
\frac{1}{n} \sum_{k=0}^{n-1} \un_A(X_k)
\stackrel{\mathit{a.s.}}{\longrightarrow}\int_A \pi \,\rmd
\mu;
\]
\item for any $a \in(0,1)$ and any continuous function $f$ in
$L_{W^a}$,
\[
\frac{1}{n} \sum_{k=0}^{n-1} f(X_k) \stackrel{\mathit{a.s.}}{\longrightarrow}
\int f   \pi  \,\rmd\mu.
\]
\end{longlist}}
\end{theo}
\begin{pf}
We check conditions \ref{coroLNGgalEEitem1},
\ref{coroLNGgalEEitem2} and \ref{coroLNGgalEEitem3} of
Corollary \ref{coroLNGgalEE} with $V \eqdef W^{a}$ for $a \in
(0,1)$, and
$\alpha\eqdef1/a$. Assumption \ref{Ageometric-ergodicity} holds and
$\limsup_n L_{\theta_n} < +\infty$ $\PP$-a.s. [see
Proposition \ref{propkey-iMCMC} and
Corollary \ref{corokey-iMCMC}\ref{corokey-iMCMCitem2}]. The drift
condition~(\ref{eqEEdriftineq}) implies that
%
\begin{equation}
\label{eqEEControlPi}
\limsup_n \pi_{\theta_n}(W) \leq\frac{\tilde b }{1 - \tilde
\lambda} \limsup_n
\theta_n(W) .\vadjust{\goodbreak}
\end{equation}
Since $W$ is continuous, the assumption \ref{EProcY}(a)
implies that $\limsup_n \theta_n(W) < \infty$ $\PP$-a.s. Hence, condition
\ref{coroLNGgalEEitem1} of Corollary \ref{coroLNGgalEE} holds.
Corollary \ref{corokey-iMCMC}\ref{corokey-iMCMCitem1} implies the
condition \ref{coroLNGgalEEitem2} of Corollary \ref{coroLNGgalEE}.
The definition (\ref{eqdefiDVnorm}) of $D_V$ implies
\[
D_V(\theta_k, \theta_{k-1}) \leq2 \upsilon\| \theta_k - \theta
_{k-1} \|_{V} \leq
\frac{2}{k+1} \theta_{k-1}(V) + \frac{2}{k+1} V(Y_{k}) .
\]
Hence, under \ref{EProcY}(a), condition
\ref{coroLNGgalEEitem3} of Corollary \ref{coroLNGgalEE} holds if
$\sum_k k^{-2} V(X_k) < +\infty$ and $ \sum_k k^{-2} V(X_k) V(Y_{k})
< +\infty$
$\PP$-a.s. The first series converges since, by
Corollary \ref{corokey-iMCMC}\ref{corokey-iMCMCitem1}, $\sup_k
\PE[V(X_k)] < +\infty$. For the second
series, it is sufficient to prove that $\sum_k k^{-2/p} V^{1/p}(X_k)
V^{1/p}(Y_{k}) < +\infty$ w.p.1 with $p \eqdef(2a ) \vee1$. We have by
the Cauchy--Schwarz inequality
\begin{eqnarray*}
\PE[ V^{1/p}(Y_{k}) V^{1/p}(X_k) ] &\leq&
\PE[V^{2/p}(Y_{k})
]^{1/2} \PE[ V^{2/p}(X_k) ]^{1/2} \\
&\leq&\PE[V^{1/a}(Y_{k}) ]^{1/2} \PE[ V^{1/a}(X_k)
]^{1/2}\\
&=& \PE[W(Y_{k})]^{1/2}  \PE[W(X_k)]^{1/2} .
\end{eqnarray*}
The RHS is finite under \ref{EProcY}(b) and
Corollary \ref{corokey-iMCMC}\ref{corokey-iMCMCitem1}. Then, this
concludes the proof of condition \ref{coroLNGgalEEitem3} of
Corollary \ref{coroLNGgalEE}.

It remains to prove that $\lim_n \pi_{\theta_n}(f) = \pi(f)$ $\PP
$-a.s. By
Proposition \ref{propITaspi}, this property holds for any bounded continuous
function $f$ and any set~$A$ such that $\int_{\partial A} \pi  \,\rmd
\mu=0$.
We proved that there exists $\alpha>1$ such that\break $\limsup_n
\pi_{\theta_n}(V^{\alpha}) + \pi(V^\alpha) < +\infty$ [see
(\ref{eqEEControlPi})]. Classical truncation arguments imply that
$\lim_n
\pi_{\theta_n}(f)$ exists $\PP$-a.s. for any continuous function $f
\in\mathcal{L}_V$ [see, e.g., \citet{billingsley1999},
Theorem 3.5,
or similar arguments in the proof of
Proposition \ref{propAscoli}].
\end{pf}

To summarize the above discussions, the process $\Xproc$ has uniformly bounded
$W$-moments (see Corollary \ref{corokey-iMCMC}), the distribution of $X_n$
converges to $\pi$ as $n \to+\infty$ (Theorem \ref
{theoEECvgMarginal}) and a
strong law of large numbers is satisfied for a wide family of functions
(Theorem \ref{theoEECvgResults}). The results are obtained provided the
auxiliary process also possesses uniformly bounded $W$-moments and
satisfies a
strong law of large numbers (see \ref{EProcY}). Repeated
applications of
this result provides sufficient conditions for the interacting
tempering with
multiple stages to be ergodic and to satisfy a strong law of large numbers.
Recall that IT algorithm defines recursively $K$ random sequences
$X^{(i)} =
\{X^{i}_n, n \geq0 \}$ for $i \in\{1, \ldots, K \}$ such that $X^{(i)}$
targets the distribution proportional to $\pi^{1/T_i}$. We are
interested in
$X^{(1)}$ which targets $\pi^{1/T_1} =\pi$. The proof of
Theorem \ref{theoEEKstage} is in Section \ref{secprooftheoEEKstage}.
%
\begin{theo}
\label{theoEEKstage} Let $(\Xset, \Xsigma)$ be a Polish space, and
$\pi$ be a density (w.r.t. a~measure $\mu$) satisfying \ref{EESPi}.
Choose $T_\star>1$ and $T_1=1 < T_2 < \cdots< T_K < T_\star$. Assume
that for any $i \in\{1, \ldots, K-1 \}$, there exists a $\pi$-irreducible
Feller transition kernel $P^{(i)}$ on $(\Xset, \Xsigma)$ such that:
{\renewcommand\thelonglist{(\roman{longlist})}
\renewcommand\labellonglist{\thelonglist}
\begin{longlist}
\item\label{itemEEKlevel5} $\pi^{1/T_i} P^{(i)} =
\pi^{1/T_i}$,\vadjust{\goodbreak}
\item\label{itemEEKlevel4} for any $s \in(0, 1/T_{i} )$, there
exist $\lambda^{(i)} \in(0,1)$
and $b^{(i)} < +\infty$ such that $P^{(i)} U_s \leq\lambda^{(i)} U_s +
b^{(i)}$ where $U_s \propto\pi^{- s}$.

\noindent Assume in addition that there exists $\tilde T \in(T_K, T_\star)$
such that:
\item\label{itemEEKlevel3} $\int\pi^{1 /T_K - 1/\tilde T} \,\rmd
\mu<
+\infty$,
\item\label{itemEEKlevel1} for any continuous function in
$\mathcal{L}_{\pi^{-1/
\tilde T}}$,
\[
n^{-1} \sum_{k=1}^n f\bigl(X_k^{(K)}\bigr) \stackrel{\mathit{a.s.}}{\longrightarrow}\int f \frac{\pi
^{1/T_K}}{\int
\pi^{1/T_K} \,\rmd\mu} \,\rmd\mu,
\]
\item\label{itemEEKlevel2} $\sup_n \PE[ \pi^{-1/\tilde
T}(X_n^{(K)})] < \infty$.
\end{longlist}}

Finally, assume that for any $i \in\{1,\ldots, K-1\}$,
$\PE[\pi^{-1/\tilde T}(X_0^{(i)}) ]< +\infty$. Then,
for any
continuous function $f$ in $\mathcal{L}_{\pi^{-1/T_\star}}$,
\[
n^{-1} \sum_{k=1}^n f\bigl(X_k^{(1)}\bigr) \stackrel{\mathit{a.s.}}{\longrightarrow}\int f  \pi \,\rmd\mu.
\]
\end{theo}

Note that since convergence holds for any continuous function $f$ in
$\mathcal{L}_{\pi^{-1/T_\star}}$, it also holds with $f = \un_A$
where $A$ is a
measurable set such that $\int_{\partial A} \pi \,\rmd\mu=0$.

We conclude this section by an example of SRWM-based interacting tempering
algorithm, for which the conditions of Theorem \ref{theoEEKstage}
hold. The
proof is in Section~\ref{secproofpropEEKstageSRWM}.
%
\begin{prop}
\label{propEEKstageSRWM}
Let $\pi$ be a super-exponential density on $\Xset= \Rset^d$ with regular
contours (i.e., satisfying \ref{AMPi}). Let $T_\star\in(1,
+\infty)$ and
choose a temperature ladder $1= T_1 < \cdots< T_K < T_\star$.
Consider the
$K$-stages interacting tempering algorithm with:
\begin{itemize}
\item for $i \in\{1, \ldots, K-1 \}$, $P^{(i)}$ is a SRWM transition
kernel with invariant distribution proportional to $\pi^{1/T_i}$ and
proposal distribution $\mathcal{N}_d(0, \Sigma^{(i)})$,
\item$\{X_n^{(K)}, n \geq0 \}$ is a SRWM Markov chain with invariant
distribution proportional to $\pi^{1/T_K}$ and proposal distribution
$\mathcal{N}_d(0, \Sigma^{(K)})$.
\end{itemize}
Finally, assume that for any $i \in\{1, \ldots, K \}$,
$\PE[\pi^{-1/T_\star}(X_0^{(i)}) ] < +\infty$. Then,~for any continuous function $f
\in\mathcal{L}_{\pi^{-1/T_\star}}$, $n^{-1} \sum_{k=1}^n
f(X_k^{(1)}) \stackrel{\mathit{a.s.}}{\longrightarrow}\pi(f)$ as \mbox{$n
\to+\infty$}.
\end{prop}

\section{\texorpdfstring{Proofs of Section \protect\ref{secMainResults}}{Proofs of Section 2}}
\label{secProofs}
\subsection{\texorpdfstring{Proof of Theorem \protect\ref{theoCvgMarg}}{Proof of Theorem 2.1}}
\label{secprooftheoCvgMarg}
We preface the proof by a lemma, which is proved in
\citet{atchadefortmoulinespriouret2009}, Proposition
1.7.1.
%
\begin{lemma}\label{lemDiminishAdap}
For any integers $n,N>0$,
\[
\sup_{\| f \|_{\infty} \leq1 }\bigl|\mathbb{E}[f(X_{n+N}) | \F _n
]-P_{\theta
_n}^Nf(X_n)\bigr|\leq\sum_{j=1}^{N-1}
\mathbb{E}[D(\theta_{n+j},\theta_n) | \F_n ] ,\qquad  \PP\mbox{-a.s.}\vadjust{\goodbreak}
\]
\end{lemma}
%
%
\begin{pf*}{Proof of Theorem \ref{theoCvgMarg}}
Let $f$ be a bounded nonnegative function.
Without loss of generality, assume that $\| f \|_{\infty} \leq1$. For any
$N \leq
n$,
%
\begin{eqnarray}\label{eqprooftheoCvgMargTool1}\qquad
| \PE[f(X_n) ] - \alpha| &\leq&|
\PE[f(X_n) - P_{\theta_{n-N}}^N f(X_{n-N}) ] |
+ | \PE[ \pi_{\theta_{n-N}}(f) - \alpha]|\nonumber\\[-10pt]\\[-10pt]
&&{} + |
\PE[ P_{\theta_{n-N}}^N f(X_{n-N}) - \pi_{\theta_{n-N}}(f)
]|.
\nonumber
\end{eqnarray}
Let $\varepsilon>0$. By setting $N =r_\varepsilon(n)$ where the
sequence $\{
r_\varepsilon(n) , n \geq0\}$ is as in~\ref{AA2}(a), the
third term on the RHS in (\ref{eqprooftheoCvgMargTool1}) is
bounded by
\[
\PE\bigl[ \bigl\|P_{\theta_{n-r_\varepsilon(n)}}^{r_\varepsilon(n)}
\bigl(X_{n-r_\varepsilon(n)}, \cdot\bigr) - \pi_{\theta_{n-r_\varepsilon
(n)}}\bigr\|_{\mathrm{TV}} \bigr]
.
\]
Under \ref{AA2}(a), for any large $n$ this
expectation is
upper bounded by $ \varepsilon$. Lemma~\ref{lemDiminishAdap} shows that
\[
\bigl| \PE\bigl[f(X_n) - P_{\theta_{n-r\varepsilon
(n)}}^{r_\varepsilon(n)}
f\bigl(X_{n-r_\varepsilon(n)}\bigr) \bigr] \bigr| \leq\sum
_{j=1}^{r_\varepsilon(n)-1}
\PE\bigl[ D\bigl(\theta_{n-r_\varepsilon(n)+j}, \theta_{n-r_\varepsilon
(n)}\bigr) \bigr] .
\]
Under \ref{AA2}(b), the RHS tends to zero
as $n \to
+\infty$. Finally, the remaining term in (\ref{eqprooftheoCvgMargTool1})
converges to zero, as a consequence of the a.s. convergence of
$\{\pi_{\theta_n}(f), n\geq0\}$ to $\alpha$, and of the property
$\lim_n n
-r_\varepsilon(n) = +\infty$.
\end{pf*}

\subsection{\texorpdfstring{Proof of Lemma \protect\ref{lemAMStability}}{Proof of Lemma 2.5}}
\label{secProofAM}
The proof of \ref{propAMBoundW} follows by iterating the drift inequality
in \citet{saksmanvihola2010}, Proposition 15. We now prove~\ref{propAMStabilityW}. Saksman and Vihola
[(\citeyear{saksmanvihola2010}), Proposition 15]
implies that there exists a~constant $c$ such that on the set $\{\sup
_{k \leq
n-1} k^{-\tau} | \theta_k | \leq t \}$,
\[
\sup_{k \leq n-1}  \lambda_{\theta_k} \leq1 - ( c t^{d
\gamma/2}
k^{\tau d \gamma/2} )^{-1} \leq1 - ( c t^{d \gamma/2}
n^{\tau d \gamma/2} )^{-1} ,\qquad \PP\mbox{-a.s.}
\]
Then by iterating the drift inequality in Saksman and Vihola
[(\citeyear{saksmanvihola2010}), Proposition~15] this yields
\begin{eqnarray*}
&&
\PE\bigl[W(X_n) \un_{\sup_{k \leq n-1} k^{-\tau} | \theta_k |
\leq t
}\bigr] \\
&&\qquad\leq\PE[W(X_0) ] + b \sum_{k=0}^{n-1} \bigl( 1 -
( c t^{d
\gamma/2}
n^{\tau d \gamma/2} )^{-1} \bigr)^k \\
&&\qquad\leq\PE[W(X_0) ] + b ( c t^{d \gamma/2} n^{\tau d
\gamma/2} ) .
\end{eqnarray*}
The last assertion follows from (\ref{eqAMStatbilitytheta}),
\ref{propAMStabilityW}, and the Markov inequality: let~$\varepsilon,\allowbreak\tau>0$;
choose $t_\varepsilon$ and $\tau'>0$ such that $\tau-\tau' d \gamma
/2 >0$ and
$\PP({\sup_{n \geq1}} |\theta_n| n^{-\tau'} \geq t_\varepsilon)
\leq\varepsilon/2$. Then
\begin{eqnarray*}
&&\PP\Bigl[\sup_n n^{-1-\tau} W(X_n) \geq M \Bigr] \\
&&\qquad\leq\varepsilon/2 + \PP\Bigl[\sup_n
n^{-1-\tau} W(X_n) \geq M, \sup_{n \geq1} |\theta_n| n^{-\tau'}
\leq t_\varepsilon\Bigr] \\
&&\qquad\leq\varepsilon/2 + \frac{1}{M} \PE\Bigl[\sup_{n \geq1}
n^{-1-\tau} W(X_n)
\un_{\sup_{n \geq1} |\theta_n| n^{-\tau'} \leq t_\varepsilon}
\Bigr] \\
&&\qquad \leq\varepsilon/2 + \frac{C}{M} \sum_{n \geq1}
\frac{1}{n^{1 +\tau}} n^{\tau' d \gamma/2}
\end{eqnarray*}
for some constant $C$, and the RHS is upper bounded by $\varepsilon$
for large
enough~$M$.\vspace*{-2pt}

\subsection{\texorpdfstring{Proof of Theorem \protect\ref{theoLGNgal}}{Proof of Theorem 2.7}}
\label{secprooftheoLGNgal}
The proof of Theorem \ref{theoLGNgal} is prefaced by lemmas on the
regularity in $\theta$ of the invariant distribution $\pi_\theta$
and on the function
$\hF_\theta$ solution of the Poisson equation $ \hF_\theta-
P_\theta\hF_\theta= F(\cdot,\theta) -
\pi_\theta(F(\cdot,\theta))$.

Under \ref{Ageometric-ergodicity}, $\hat F_\theta(x) \eqdef\sum_n
P_\theta^n
\{F(\cdot,\theta) - \pi_\theta(F(\cdot,\theta)) \}(x)$ exists for
all $x \in\Xset$, solves
the Poisson equation, and by Lemma \ref{lemBoundCandRho}
%
\begin{equation}
\label{eqUpperBoundPoisson}
| \hat F_\theta(x) | \leq\| F(\cdot, \theta) \|_{V}
L_\theta^2   V(x) ,
\end{equation}
where $L_\theta$ is defined in (\ref{eqDefinitionLtheta}).

The following lemma is adapted from \citet{andrieujasradoucetdelmoral2008}.
A detailed proof is given in Section 3 of the
supplemental paper [\citet{supplement}].\vspace*{-2pt}
%
\begin{lemma}
\label{lemReguThetaPoisson} Assume \ref{Ageometric-ergodicity}. For any
$\theta\in\Theta$, let $F_\theta\dvtx \Xset\to \Rset^+$ be a measurable
function such that $\sup_\theta\| F_\theta \|_{V} < +\infty$ and define
$\hF_\theta\stackrel{\mathit{def}}{=} \sum_{n \geq0}
P_\theta^n\{F_\theta- \pi_\theta(F_\theta) \}$. For any $\theta,
\theta' \in\Theta$,
\[
\| \pi_\theta- \pi_{\theta'} \|_{V} \leq L_{\theta'}^2 \{
\pi_\theta(V) + L_\theta^2  V(x)
\}  D_{V}(\theta,\theta')
\]
and
\begin{eqnarray*}
| P_\theta\hF_\theta- P_{\theta'} \hF_{\theta'} |_V
&\leq&{\sup_{\theta\in\Tset}}
\| F_\theta\|_{V}  L_{\theta'}^2  \bigl( L_{\theta}
D_V(\theta,\theta') + \| \pi_\theta- \pi_{\theta'} \|_{V}
\bigr) \\[-3pt]
&&{} + L_{\theta'}^2  \| F_\theta- F_{\theta'} \|_{V} ,
\end{eqnarray*}
where $L_\theta$ is given by (\ref{eqDefinitionLtheta}).\vspace*{-2pt}
\end{lemma}
\begin{pf*}{Proof of Theorem \ref{theoLGNgal}}
We denote by $L$ the limit $\lim_{n } \int\pi_{\theta_n}(\rmd x)
F(\theta_n,\allowbreak x)$. We write $ \frac{1}{n} \sum_{k=0}^{n-1} F(X_k,
\theta_k) - L =
\sum_{i=1}^4 T_{i,n}$ with
\begin{eqnarray*}
T_{1,n} &\eqdef& \frac{1}{n}F(X_0, \theta_0) - \frac{L}{n}, \\
T_{2,n} &\eqdef& \frac{1}{n} \sum_{k=1}^{n-1} \{F(X_k, \theta_k)
-F(X_k, \theta_{k-1}) \}, \\
T_{3,n} &\eqdef& \frac{1}{n}\sum_{k=1}^{n-1} \biggl\{F(X_k,
\theta_{k-1}) - \int\pi_{\theta_{k-1}}(\rmd x) F(x,\theta_{k-1})
\biggr\}, \\
T_{4,n} &\eqdef& \frac{1}{n}\sum_{k=0}^{n-2} \biggl\{ \int\pi
_{\theta_{k}}(\rmd x)
F(x,\theta_{k}) - L \biggr\}.
\end{eqnarray*}
Consider first $T_{1,n}$. Since $|F(X_0,\theta_0)| < +\infty$ $\PP
$-a.s., $\lim_{n
\to\infty} T_{1,n} =0 $ $\PP$-a.s. Under conditions
\ref{theoLGNgalitem2} [resp., \ref{theoLGNgalitem3}],
$T_{2,n}$ (resp.,
$T_{4,n}$) converges to zero a.s. (for $T_{2,n}$, note that $L_\theta
\geq1$ by
definition). Consider finally $T_{3,n}$:
\[
\frac{1}{n} \sum_{k=1}^{n-1} \biggl\{ F(X_k, \theta_{k-1}) - \int
\pi_{\theta_{k-1}}(\rmd x) F(x,\theta_{k-1}) \biggr\} = M_n + R_n +
\tilde R_n
\]
with $\hF_\theta(x) \eqdef\sum_{n \geq0} P_\theta^n \{F(\cdot
,\theta) -
\pi_\theta(F(\cdot,\theta)) \}(x)$ and
\begin{eqnarray*}
M_n & \eqdef &\frac{1}{n} \sum_{k=1}^{n-1} \{ \hF_{\theta_{k-1}}
(X_k) -
P_{\theta_{k-1}} \hF_{\theta_{k-1}}(X_{k-1})\} ,\\
R_n & \eqdef &\frac{1}{n} \sum_{k=1}^{n-1} \{ P_{\theta_{k}}
\hF_{\theta_{k}}(X_{k}) -P_{\theta_{k-1}}
\hF_{\theta_{k-1}}(X_{k}) \} , \\
\tilde R_n & \eqdef &\frac{1}{n} P_{\theta_0} \hF_{\theta_0}(X_0) -
\frac{1}{n} P_{\theta_{n-1}} \hF_{\theta_{n-1}}(X_{n-1}) .
\end{eqnarray*}
By construction, $\{ \hF_{\theta_{k-1}} (X_k) - P_{\theta_{k-1}}
\hF_{\theta_{k-1}}(X_{k-1}), k \geq1 \}$ is a martingale-increment sequence.
Therefore, by Hall and Heyde [(\citeyear{hallheyde1980}), Theorem 2.18], $M_n \aslim0$ provided
that
%
\begin{equation}
\label{eqtheoLGNgaltool1}\quad
\sum_{k \geq1} \frac{1}{k^{\alpha}} \mathbb{E}[|\hF
_{\theta _{k-1}} (X_k) - P_{\theta_{k-1}} \hF_{\theta
_{k-1}}(X_{k-1}) |^\alpha | \F _{k-1} ] < +\infty,\qquad \PP
\mbox{-a.s.}
\end{equation}
Equation (\ref{eqUpperBoundPoisson}) and Jensen's inequality imply that
($\alpha
>1$)
\begin{eqnarray*}
&& \mathbb{E}[|\hF_{\theta_{k-1}} (X_k) - P_{\theta_{k-1}}
\hF _{\theta_{k-1}}(X_{k-1}) |^\alpha | \F_{k-1} ] \\
&&\qquad  \leq2^{\alpha-1} \mathbb{E}[|\hF_{\theta_{k-1}}
(X_k) |^\alpha + |P_{\theta_{k-1}}\hF_{\theta_{k-1}}
(X_{k-1}) |^\alpha | \F_{k-1} ] \\
&&\qquad  \leq2^\alpha\Bigl( \sup_\theta\| F(\cdot,\theta) \|_{V}
 L_{\theta_{k-1}}^2
\Bigr)^\alpha P_{\theta_{k-1}}V^\alpha(X_{k-1}) .
\end{eqnarray*}
Under item \ref{theoLGNgalitem1} and
\ref{AUniformBoundLGN}(b), the series
is finite
$\PP$-a.s. and this concludes the proof of~(\ref{eqtheoLGNgaltool1}).
Consider now the remainder term $R_n$. By Lemma \ref{lemReguThetaPoisson},
\begin{eqnarray*}
|R_n| &\leq&\frac{\sup_\theta\| F(\cdot,\theta) \|_{V}}{n} \\
&&{}\times\sum
_{k=1}^n L_{\theta_k}^2
L_{\theta_{k-1}}^2 \{ 1 +
\pi_{\theta_k}(V) + L_{\theta_k}^2 \}  D_V(\theta_k, \theta
_{k-1}) V(X_k) \\
&&{} + \frac{1}{n} \sum_{k=1}^n L_{\theta_{k}}^2 \| F(\cdot,\theta_k)
- F(\cdot,\theta_{k-1}) \|_{V} V(X_k) .
\end{eqnarray*}
Assumptions\vspace*{2pt} \ref{AStrongContainmentLGN},
\ref{AUniformBoundLGN}(a) and
items \ref{theoLGNgalitem1}, \ref{theoLGNgalitem2} imply that $R_n
\aslim0$. Consider finally~$\tilde R_n$. By (\ref{eqUpperBoundPoisson}),
\begin{eqnarray*}
&&\frac{1}{n} | P_{\theta_0} \hF_{\theta_0}(X_0) - P_{\theta_{n-1}}
\hF_{\theta_{n-1}}(X_{n-1}) | \\
&&\qquad\leq\frac{\sup_\theta\| F(\cdot,\theta) \|_{V}}{n}  \bigl(
L_{\theta_0}^2  P_{\theta_0}
V(X_0)+
L_{\theta_{n-1}}^2  P_{\theta_{n-1}} V(X_{n-1})\bigr) \\
&&\qquad\leq\frac{\sup_\theta\| F(\cdot,\theta) \|_{V}}{n}  \bigl(
L_{\theta_0}^2  \{ V(X_0)
+ b_{\theta_0} \}+ L_{\theta_{n-1}}^2  P_{\theta_{n-1}}
V(X_{n-1})\bigr) .
\end{eqnarray*}
Assumption \ref{AUniformBoundLGN}(b),
item \ref{theoLGNgalitem1} and the condition $ V(X_0) < +\infty
$ $\PP$-a.s.
imply that $\tilde R_n \aslim0$.
\end{pf*}

\subsection{\texorpdfstring{Proof of Theorem \protect\ref{theoaspi}}{Proof of Theorem 2.11}}
\label{secproofaspi}
We preface the proof of this theorem by a~proposition and a lemma. The
proof of
Proposition \ref{propAscoli} is postponed
to \citet{supplement}, Section 4.
%
\begin{prop}
\label{propAscoli}
Let $\Xset$ be a Polish space endowed with its Borel $\sigma$-field
$\Xsigma$.
Let $\mu$ and $\{\mu_n,n \geq1 \}$ be probability distributions on
$(\Xset,
\Xsigma)$. Let $\{h_n, n \geq0\}$ be an equicontinuous family of functions
from $\Xset$ to $\Rset$. Assume:
{\renewcommand\thelonglist{(\roman{longlist})}
\renewcommand\labellonglist{\thelonglist}
\begin{longlist}
\item the sequence $\{\mu_n,n \geq0 \}$ converges weakly to $\mu$,
\item for any $x \in\Xset$, $\lim_{n \to\infty} h_n(x)$ exists,
and there
exists $\alpha>1$ such that $\sup_n \mu_n(|h_n|^\alpha) + \mu
(|\lim_{n}
h_n|) < +\infty$.
\end{longlist}}

\noindent Then, $\mu_n(h_n) \to\mu( \lim_{n } h_n)$.
\end{prop}
%
\begin{lemma}
\label{lemEECvgIteratedKernels}
Let $\Xset$ be a Polish space endowed with its Borel $\sigma$-field
$\Xsigma$.
Let $\{P_\theta, \theta\in\Theta\}$ be a family of transition
kernels on $(\Xset,
\Xsigma)$ and $\{\theta_n, n\geq0\}$ be a $\Theta$-valued random
sequence on
$(\Omega, \mathcal{A}, \PP)$. Assume conditions
\ref{theoaspiitem3} and \ref{theoaspiitem2} of Theorem~\ref
{theoaspi}.
Then, there exists a $\PP$-full set $\Omega_\star$ such that for any
$\omega
\in\Omega_\star$, $x \in\Xset$ and $k \geq1$, the probability
distributions
$\{P_{\theta_n(\omega)}^k(x,\cdot), n\geq0\}$ converge weakly to
$P_{\theta_\star}^k(x,\cdot)$.
\end{lemma}
\begin{pf}
We prove, by induction on $k$, that there exists a $\PP$-full set
$\Omega_k$
such that for any $\omega\in\Omega_k$ and $x \in\Xset$, the probability
distributions $\{P_{\theta_n(\omega)}^k(x,\cdot),\allowbreak n\geq0\}$
converge weakly to
$P_{\theta_\star}^k(x,\cdot)$. The proof is then concluded by setting
$\Omega_\star= \bigcap_k \Omega_k$.

Consider the case $k=1$. By condition \ref{theoaspiitem2} of Theorem
\ref{theoaspi}, for any $x \in\Xset$ there exists a $\PP $-full set
$\Omega_x$ such that for any $\omega\in\Omega_x$,
$\{P_{\theta_n(\omega)}(x,\cdot), n\geq0\}$ converges weakly to
$P_{\theta_\star} (x,\cdot)$. Since $\Xset$ is Polish, it admits a
countable dense subset~$\mathcal{D}$. Therefore, there exists a
$\PP$-full set $\Omega_{\mathcal{D}}$ such that for any
$\omega\in\Omega _{\mathcal{D}}$ and any $x \in\mathcal{D}$, $\{
P_{\theta_n(\omega)}(x,\cdot), n\geq0\}$ converges weakly to
$P_{\theta_\star}(x, \cdot)$. Under condition \ref{theoaspiitem3} of
Theorem \ref{theoaspi}, for any bounded continuous function $f$, the
family of functions $\{\bar P_\theta f \eqdef P_\theta f -
P_{\theta_\star} f, \theta\in\Theta\}$ is equicontinuous. For any
$\varepsilon>0$ and any $x \in\Xset$, there thus exists
$x_\varepsilon\in\mathcal{D}$ such that for any $\theta\in \Theta$,
$|\bar P_{\theta} f(x) - \bar P_{\theta} f(x_\varepsilon)|
\leq\varepsilon $. Hence, for any $\omega\in\Omega_{\mathcal{D}}$ and
any bounded continuous function $f$,
\begin{eqnarray*}
\bigl| \bar P_{\theta_n(\omega)} f(x) \bigr| &\leq&\bigl| \bar
P_{\theta_n(\omega)}
f(x_\varepsilon) \bigr| + \bigl|\bar
P_{\theta_n(\omega)} f(x) - \bar P_{\theta_n(\omega)}
f(x_\varepsilon)\bigr| \\
&\leq&\bigl| P_{\theta_n(\omega)} f(x_\varepsilon) - P_{\theta
_\star} f(x_\varepsilon) \bigr| + \varepsilon.
\end{eqnarray*}
This implies that $\limsup_n | \bar P_{\theta_n(\omega)} f(x)
| \leq
\varepsilon$. Since $\varepsilon$ was arbitrary, it follows
$\{P_{\theta_n(\omega)}(x,\cdot), n\geq0\}$ converges weakly to
$P_{\theta_\star}(x,\cdot)$ for any $x$. Hence, we set $\Omega_1 =
\Omega_{\mathcal{D}}$.

Assume that the property holds for $k \geq1$. We write for any bounded and
continuous function $f$
%
\begin{eqnarray}\label{eqCvgIteratedKernelsTool1}\quad
P_{\theta_n(\omega)}^{k+1} f(x) - P_{\theta_\star}^{k+1} f(x) &=&
\int
\bigl(P_{\theta_n(\omega)}^{k}(x,\rmd y ) - P_{\theta_\star
}^{k}(x,\rmd y ) \bigr)
P_{\theta_\star}f(y) \nonumber\\[-8pt]\\[-8pt]
&&{} + \int P_{\theta_n(\omega)}^{k}(x,\rmd y ) \bigl(P_{\theta_n(\omega
)}f(y) - P_{\theta_\star}f(y) \bigr).
\nonumber
\end{eqnarray}
By the induction assumption, there exists a $\PP$-full set $\Omega_k$
such that
for any $\omega\in\Omega_k$, $x \in\Xset$ and any bounded continuous
function\vspace*{1pt} $h$, $ \lim_{n \to\infty} P_{\theta_n(\omega)}^k h(x) =
P_{\theta_\star}^k
h(x)$. Applied with $h = P_{\theta_\star} f$, which is continuous
under the
assumption~\ref{theoaspiitem3}, this proves that for any $\omega
\in
\Omega_k$, the first term on the RHS of (\ref{eqCvgIteratedKernelsTool1})
goes to zero. For the second term, we use Proposition \ref
{propAscoli}. Let
$\omega\in\Omega_k \cap\Omega_1$. For any $x \in\Xset$, $\{
P_{\theta_n(\omega)}^k(x,\cdot), n\geq0\}$ converges weakly to
$P_{\theta_\star}^k(x,\cdot)$. Furthermore, the family of bounded
functions $\{
P_{\theta_n(\omega)}f - P_{\theta_\star} f, n\geq0\}$ is
equicontinuous and, since
$\omega\in\Omega_1$, $\lim_{n \to\infty} P_{\theta_n(\omega
)}f(y) - P_{\theta_\star}
f(y) =0$ for any $y \in\Xset$. Proposition~\ref{propAscoli} thus implies
that the second term on the RHS of (\ref{eqCvgIteratedKernelsTool1})
converges to zero, for any bounded continuous function $f$. The above
discussion proves that $\Omega_{k+1} = \Omega_k \cap\Omega_1 =
\Omega_1$, and
concludes the induction.
\end{pf}
\begin{pf*}{Proof of Theorem \ref{theoaspi}}
Fix $x \in\Xset$. Let $f$ be a bounded continuous function on $\Xset$.
Under \ref{Ageometric-ergodicity}, we have by
Lemma \ref{lemBoundCandRho}
\begin{eqnarray*}
&&
{\limsup_n} | \pi_{\theta_n}(f) - P_{\theta_n}^k f(x) +
P_{\theta_\star}^k f(x)-
\pi_{\theta_\star}(f) | \\
&&\qquad
\leq\Bigl( \limsup_n C_{\theta_n}  [\limsup_n \rho_{\theta
_n}]^k + C_{\theta_\star}
\rho_{\theta_\star}^k \Bigr) V(x) .
\end{eqnarray*}
By Lemma \ref{lemBoundCandRho} and condition \ref{theoaspiitem1},
$\limsup_n C_{\theta_n} <+\infty$ and $\limsup_n \rho_{\theta_n}
<1 $ $\PP\mbox{-a.s.}$; then,
there exists a $\PP$-full set $\Omega_\star''$ such that for any
$\omega\in
\Omega_\star''$, there exists~$k(\omega)$ such that
\[
\limsup_n \bigl| \pi_{\theta_n(\omega)}(f) - P_{\theta_n(\omega
)}^{k(\omega)} f(x) +
P_{\theta_\star}^{k(\omega)} f(x)- \pi_{\theta_\star}(f) \bigr|
\leq\varepsilon.
\]
Note\vspace*{1pt} that $\Omega_\star''$ does not depend upon $x$ and $f$. By Lemma
\ref{lemEECvgIteratedKernels}, there exists a~$\PP$-full set
$\Omega_\star$ such that $\lim_{n \to\infty} P_{\theta_n(\omega )}^{k}
f(x) = P_{\theta_\star}^{k} f(x)$ for any\vspace*{1pt}
$\omega\in\Omega_\star$, any $x \in\Xset$, any $k \geq1$ and any
bounded continuous function $f$. The proof is concluded by setting
$\Omega_\star' = \Omega_\star'' \cap\Omega_\star$.
\end{pf*}

\section{\texorpdfstring{Proofs of Section \protect\ref{secEE}}{Proofs of Section 3}}
\label{secProofsEE}

\subsection{\texorpdfstring{Proof of Corollary \protect\ref{corokey-iMCMC}}{Proof of Corollary 3.2}}
\label{secproofcorokey-iMCMC}

\ref{corokey-iMCMCitem1} By iterating the drift inequality (\ref
{eqEEdriftineq}), we obtain
\[
\PE[W(X_n) ] \leq\tilde\lambda^n \PE[W(X_0)
] + \tilde
b \sum_{k=0}^{n-1} \tilde\lambda^k \PE[ \theta
_{n-k}(W)] .
\]
Under \ref{EProcY}(b), $\sup_{k \geq0} \PE[
\theta_{k}(W)] < +\infty$ so that
%
\begin{equation}
\label{eqIMCMC-key2}
\PE[W(X_n) ] \leq\tilde\lambda^n \PE[W(X_0)
] +
\frac{\tilde b}{1 - \tilde\lambda} \sup_{k \geq0} \PE[
\theta_{k}(W)]
.
\end{equation}

\ref{corokey-iMCMCitem2} Since $W$ is a continuous function,
\ref{EProcY}(a) implies that $\limsup_n \theta _n(W) < +\infty$,
$\PP$-a.s. Consequently, $\limsup_{n } L_{\theta_n} < +\infty$,
$\PP$-a.s. by Lemma \ref{lemBoundCandRho} and
Proposition~\ref{propkey-iMCMC}.

\subsection{\texorpdfstring{Proof of Proposition \protect\ref{propITaspi}}{Proof of Proposition 3.3}}
\label{secproofpropITaspi}
We check the conditions of Theorem \ref{theoaspi}.
Condition \ref{theoaspiitem1} of Theorem \ref{theoaspi} holds by
Corollary \ref{corokey-iMCMC}.

The proof of condition \ref{theoaspiitem3} of Theorem \ref
{theoaspi} is
a consequence of the following lemma.
%
\begin{lemma}
\label{lemmaFeller}
Let $f$ be a function on $\Xset$ such that $\| f \pi^\beta\|
_{\infty} <
+\infty$. For any $x,x' \in\Xset$ such that $\pi(x) >0, \pi(x')>0$,
\begin{eqnarray*}
\sup_{\theta\in\Theta} | P_\theta f(x) - P_\theta f(x')
| &\leq&| P f(x)
- P f(x') | + |f(x) - f(x') | \\
&&{} + 2 \| f   \pi ^\beta\|_{\infty}  | \pi^{-\beta}(x) -
\pi^{-\beta}(x') | .
\end{eqnarray*}
\end{lemma}

\begin{pf}
By definition of the transition kernel $P_\theta$, it is easily
checked that
%
\begin{eqnarray}\label{eqlemmaFellertool1}
&&
P_\theta f(x) - P_\theta f(x') \nonumber\\
&&\qquad= \upsilon\int\{\alpha(x,y) -
\alpha(x',y)
\} \bigl(f(y)
-f(x') \bigr) \theta(\rmd y) \\
&&\qquad\quad{}+ (1-\upsilon) \bigl( Pf(x) - Pf(x') \bigr) + \upsilon\bigl( f(x)
- f(x') \bigr) A(\theta,x),
\nonumber
\end{eqnarray}
where $A(\theta,x) \eqdef1 - \int\alpha(x,y) \theta(\rmd y)$.
Since $0 \leq\alpha(x,y) \leq1$, we have
\[
\bigl| \upsilon\bigl( f(x) - f(x') \bigr) A(\theta,x) \bigr|
\leq| f(x) - f(x') | .
\]
We can assume w.l.o.g. that $\pi(x) \leq\pi(x')$. By definition of
the ratio
$\alpha$, we have
\begin{eqnarray*}
\alpha(x,y) - \alpha(x',y) &=& \un_{\{ \pi(x) \leq\pi(y) \leq\pi
(x') \}} \bigl(\pi^{-\beta}(y)
- \pi^{-\beta}(x') \bigr) \pi^\beta(y) \\
&&{}
+\bigl(\pi^{-\beta}(x) - \pi^{-\beta}(x') \bigr) \un_{\{\pi(y)
\leq\pi(x) \leq\pi(x') \}} \pi^\beta(y) ,
\end{eqnarray*}
showing that $
|\alpha(x,y) - \alpha(x',y) | \leq(\pi^{-\beta
}(x) - \pi^{-\beta}(x') )
\pi^\beta(y) \un_{\{ \pi(y) \leq\pi(x') \}}$.
The proof is concluded by noting that
\begin{eqnarray*}
&&
\int|\alpha(x,y) - \alpha(x',y) | | f(y)
-f(x')|\theta(\rmd y) \\
&&\qquad\leq2 \Bigl( \sup_{\Xset} |f| \pi^\beta\Bigr) \bigl(\pi
^{-\beta}(x) - \pi^{-\beta}(x') \bigr) .
\end{eqnarray*}
\upqed\end{pf}

The most delicate part consists in establishing
condition \ref{theoaspiitem2} of Theorem~\ref{theoaspi}. The
proof relies
on the following result which is an extension of the Varadarajan theorem
[\citet{dudley2002}, Theorem 11.4.1]. The proof of
Proposition~\ref{propEEGeneralTensionAS} is detailed in
Section 5 of the supplemental paper
[\citet{supplement}].
%
\begin{prop}
\label{propEEGeneralTensionAS}
Let $(\Uset,d)$ be a metric space equipped with its Borel $\sigma$-field
$\mathcal{B}(\Uset)$. Let $(\Omega, \mathcal{A}, \PP)$ be a
probability space,
$\mu$ be a distribution on $(\Uset, \mathcal{B}(\Uset))$ and $\{
K_n, n\geq0\}$ be
a family of Markov transition kernels $K_n \dvtx \Omega\times\mathcal
{B}(\Uset)
\to[0,1]$. Assume that, for any $f \in\Cb(\Uset,d)$
\[
\Omega_f \stackrel{\mathit{def}}{=}\Bigl\{\omega\in\Omega\dvtx
{\limsup_{n \to\infty}} | K_n(\omega,f) - \mu(f) | =0 \Bigr\}
\]
is a $\PP$-full set. Then
\[
\Bigl\{\omega\in\Omega\dvtx \forall f \in\Cb(\Uset,d)\  \limsup
_{n \to\infty} | K_n(\omega,f) -
\mu(f) | =0 \Bigr\}
\]
is a $\PP$-full set.
\end{prop}
\begin{pf*}{Proof of \ref{theoaspiitem2} of Theorem \ref{theoaspi}}
We check the conditions of Proposition \ref{propEEGeneralTensionAS}
with $\mu_n =
P_{\theta_n}(x,\cdot)$ and $\mu= P_{\theta_\star}(x,\cdot)$. For
any $x \in\Xset$,
and $f \in\Cb(\Xset)$, $y \mapsto\alpha(x,y)$ and $y \mapsto
\alpha(x,y)
f(y)$ are continuous. Thus, \ref{EProcY}(a) implies that
$P_{\theta_n}f(x) \aslim P_{\theta_\star}f(x)$ and $\Omega_f$ is a
$\PP$-full set.
\end{pf*}

\subsection{\texorpdfstring{Proof of Theorem \protect\ref{theoEEKstage}}{Proof of Theorem 3.6}}
\label{secprooftheoEEKstage}
Set $\alpha_0 = 1$ and choose $\alpha_l >1$ such that $\tilde T
\times
\prod_{l=0}^{K} \alpha_l = T_\star$. The proof is by induction on
$i$ for $i
=K$ down to $i=2$.

Set $W^{(K-1)} \eqdef\pi^{-T_\star^{-1} \prod_{l=0}^{K-1} \alpha
_l} =
\pi^{-1/(\tilde T \alpha_K)}$ and $\pi^{(K-1)}$ be the probability
distribution proportional to $\pi^{1/T_{K-1}}$. Under the stated
assumptions, Theorem \ref{theoEECvgResults} applies with $Y
\leftarrow
X^{(K)}$ and $X \leftarrow X^{(K-1)}$: for any continuous function $f$ in
$\mathcal{L}_{W^{(K-1)}}$, $n^{-1} \sum_{k=1}^n f(X_k^{(K-1)}) \aslim
\pi^{(K-1)}(f)$.

Assume\vspace*{1pt} Theorem \ref{theoEECvgResults} holds with $Y \leftarrow
X^{(i+1)}$ and $X \leftarrow X^{(i)}$ for some $i \in\{2, \ldots,
K-1\}$:
for any continuous function $f$ in $\mathcal{L}_{W^{(i)}}$, $n^{-1}
\sum_{k=1}^n
f(X_k^{(i)}) \aslim\pi^{(i)}(f)$ where $W^{(i)} \eqdef\pi^{-T_\star^{-1}
\prod_{l=0}^{i} \alpha_l}$ and $\pi^{(i)} \propto\pi^{1/T_i}$. We apply
the above results with
\begin{eqnarray*}
\pi&\leftarrow&\pi^{1/T_{i-1}} ,\qquad
\theta_\star\leftarrow\pi^{1/T_i},\qquad
P \leftarrow P^{(i-1)} ,\qquad\\
T &\leftarrow&\frac{T_i}{T_{i-1}} , \qquad
W \leftarrow W^{- T_\star^{-1} \prod_{l=0}^i \alpha_l} .
\end{eqnarray*}
We thus have that $n^{-1} \sum_{k=1}^n f(X_k^{(i-1)}) \aslim\pi^{(i-1)}(f)$
for any continuous function $f$ in $\mathcal{L}_{W^{(i-1)}}$, where
\[
W^{(i-1)} \eqdef\pi^{-T_\star^{-1} \prod_{l=0}^{i-1} \alpha_l} = \bigl\{W^{(i)}
\bigr\}^{1/\alpha_i},\qquad  \pi^{(i-1)} \propto\pi^{1/T_{i-1}} .
\]
This concludes the induction.

\subsection{\texorpdfstring{Proof of Proposition \protect\ref{propEEKstageSRWM}}{Proof of Proposition 3.7}}
\label{secproofpropEEKstageSRWM}
For any $i \in\{1, \ldots, K \}$, the transition kernels $P^{(i)}$ are
$\pi$-irreducible, aperiodic, and compact sets are $1$-small. In
addition, they
are Feller (the proof is on the same lines as the proof of
Lemma~\ref{lemmaFeller}). By Saksman and Vihola
[(\citeyear{saksmanvihola2010}), Proposition 15]
conditions \ref{itemEEKlevel5} and~\ref{itemEEKlevel4} of
Theorem \ref{theoEEKstage} are satisfied for $i \in\{1, \ldots, K
\}$. Note
that the proof of Proposition 15 in \citet{saksmanvihola2010} is in the
case $s
T_i = 1/2$ but it can be easily adapted for any $s T_i\in(0,1)$. In
the case
$i=K$, this implies that there exist $\lambda\in(0,1)$ and $b <
+\infty$ such
that
\[
P^{(K)} \tilde U \leq\lambda\tilde U + b,
\]
where $\tilde U = (\pi/ \sup_\Xset\pi)^{-1/\tilde T}$. Standard
results on
Markov chains [see, e.g., \citet{meyntweedie2009}] imply
\ref{itemEEKlevel1}. By iterating the drift inequality, we have
\[
\sup_n \PE\bigl[ \tilde U\bigl(X_n^{(K)}\bigr) \bigr] \leq\PE\bigl[
\tilde
U\bigl(X_0^{(K)}\bigr) \bigr] + \frac{b}{1-\lambda} ,
\]
thus proving \ref{itemEEKlevel2}. Finally, since $\pi$ satisfies
\ref{AMPi}, there exist positive constants~$c_i$ such that $\pi(x)
\leq c_1
\exp(-c_2 |x|)$ [see, e.g., \citet{saksmanvihola2010}, Lemma~8].
Therefore, for
any $\tau>0$, $\int\pi^{\tau}(x) \,\rmd x < +\infty$ thus showing~\ref{itemEEKlevel3}.

\section*{Acknowledgments}

The authors are indebted to the associate editor and the referees for
the numerous suggestions made to improve the manuscript.

\begin{supplement}[id=suppA]
\stitle{Supplement to paper ``Convergence of
adaptive and interacting Markov chain
Monte Carlo algorithms''}
\slink[doi]{10.1214/11-AOS938SUPP} 
\sdatatype{.pdf}
\sfilename{aos938\_supp.pdf}
\sdescription{This supplement provides a detailed proof of
Lemma \ref{lemReguThetaPoisson} and Propositions \ref{propkey-iMCMC},
\ref{propAscoli} and \ref{propEEGeneralTensionAS}. It also contains a
discussion on the setwise convergence of transition kernels.}
\end{supplement}


\printaddresses

\end{document}